\documentclass[11pt, reqno]{amsart}

\usepackage{amsmath,amssymb,amsthm,amsfonts,verbatim}
\usepackage{microtype}
\usepackage[all,2cell]{xy}
\usepackage{mathtools}
\usepackage{graphicx}
\usepackage{pinlabel}
\usepackage{hyperref}
\usepackage{mathrsfs}
\usepackage{color}
\usepackage[dvipsnames]{xcolor}
\usepackage{enumitem}
\usepackage{cite}
\usepackage{soul}

\CompileMatrices

\usepackage[top=1.2in,bottom=1.2in,left=1in,right=1in]{geometry}

\usepackage{hyperref} 
\hypersetup{          
	colorlinks=true, breaklinks, linkcolor=[RGB]{51 102 204}, filecolor=Orchid, urlcolor=[RGB]{51 102 204},
	citecolor=Orchid, linktoc=all, }
\usepackage[nameinlink]{cleveref}

\theoremstyle{plain}
\newtheorem{theorem}{Theorem}[section]
\newtheorem{maintheorem}{Theorem}
    \AddToHook{env/maintheorem/begin}{\crefalias{maintheorem}{maintheorem}}
    \crefname{maintheorem}{theorem}{theorems}

\newtheorem{maincor}[maintheorem]{Corollary}
    \AddToHook{env/maincor/begin}{\crefalias{maintheorem}{maincor}}
    \crefname{maincor}{corollary}{corollaries}

    \AddToHook{env/problem/begin}{\crefalias{theorem}{problem}}
    \crefname{problem}{problem}{problems}
\newtheorem{proposition}[theorem]{Proposition}
    \AddToHook{env/proposition/begin}{\crefalias{theorem}{proposition}}
    \crefname{proposition}{proposition}{propositions}
\newtheorem{lemma}[theorem]{Lemma}
    \AddToHook{env/lemma/begin}{\crefalias{theorem}{lemma}}
    \crefname{lemma}{lemma}{lemmas}

    \AddToHook{env/conjecture/begin}{\crefalias{theorem}{conjecture}}
    \crefname{conjecture}{conjecture}{conjectures}

    \AddToHook{env/fact/begin}{\crefalias{theorem}{fact}}
    \crefname{fact}{fact}{facts}

    \AddToHook{env/corollary/begin}{\crefalias{theorem}{corollary}}
    \crefname{corollary}{corollary}{corollaries}

    \AddToHook{env/claim/begin}{\crefalias{theorem}{claim}}
    \crefname{claim}{claim}{claims}

\theoremstyle{definition}
\newtheorem{definition}[theorem]{Definition}
    \AddToHook{env/definition/begin}{\crefalias{theorem}{definition}}
    \crefname{definition}{definition}{definitions}

    \AddToHook{env/example/begin}{\crefalias{theorem}{example}}
    \crefname{example}{example}{examples}

    \AddToHook{env/construction/begin}{\crefalias{theorem}{construction}}
    \crefname{construction}{construction}{constructions}

    \AddToHook{env/exercise/begin}{\crefalias{theorem}{exercise}}
    \crefname{exercise}{exercise}{exercises}
\newtheorem{remark}[theorem]{Remark}
    \AddToHook{env/remark/begin}{\crefalias{theorem}{remark}}
    \crefname{remark}{remark}{remarks}

    \AddToHook{env/question/begin}{\crefalias{theorem}{question}}
    \crefname{question}{question}{questions}

    \AddToHook{env/answer/begin}{\crefalias{theorem}{answer}}
    \crefname{answer}{answer}{answers}

    \AddToHook{env/convention/begin}{\crefalias{theorem}{convention}}
    \crefname{convention}{convention}{conventions}

\newcommand{\nc}{\newcommand}
\nc{\dmo}{\DeclareMathOperator}

\nc{\cA}{\mathcal{A}}
\nc{\cB}{\mathcal{B}}
\nc{\sB}{\mathscr{B}}
\nc{\BB}{\mathbb{B}}
\nc{\C}{\mathbb{C}}
\nc{\cC}{\mathcal{C}}
\nc{\sC}{\mathscr{C}}
\nc{\bd}{\mathbf{d}}
\nc{\DD}{\mathbb{D}}
\nc{\cD}{\mathcal{D}}
\nc{\sD}{\mathscr{D}}
\nc{\cE}{\mathcal{E}}
\nc{\sE}{\mathscr{E}}
\nc{\bF}{\mathbb{F}}
\nc{\cF}{\mathcal{F}}
\nc{\sF}{\mathscr{F}}
\nc{\cG}{\mathcal{G}}
\nc{\sG}{\mathscr{G}}
\nc{\cH}{\mathcal{H}}
\nc{\cI}{\mathcal{I}}
\nc{\sI}{\mathscr I}
\nc{\cK}{\mathcal{K}}
\nc{\sK}{\mathscr{K}}
\nc{\cL}{\mathcal{L}}
\nc{\sL}{\mathscr{L}}
\nc{\cM}{\mathcal{M}}
\nc{\bM}{\mathbf{M}}
\nc{\sM}{\mathscr{M}}
\nc{\N}{\mathbb{N}}
\nc{\cN}{\mathcal{N}}
\nc{\cO}{\mathcal{O}}
\nc{\bp}{\mathbf{p}}
\nc{\bP}{\mathbb{P}}
\nc{\cP}{\mathcal{P}}
\nc{\Q}{\mathbb{Q}}
\nc{\R}{\mathbb{R}}
\nc{\sR}{\mathscr{R}}
\nc{\cS}{\mathcal{S}}
\nc{\cT}{\mathcal T}
\nc{\cU}{\mathcal U}
\nc{\sU}{\mathscr U}
\nc{\bV}{\mathbb V}
\nc{\sV}{\mathscr V}
\nc{\cX}{\mathcal{X}}
\nc{\cY}{\mathcal{Y}}
\nc{\Z}{\mathbb{Z}}
\nc{\disk}{\mathbb{D}}
\nc{\hyp}{\mathbb{H}}
\renewcommand{\P}{\mathbb{P}}

\nc{\CP}{\mathbb{CP}}
\nc{\RP}{\mathbb{RP}}
\dmo{\Mod}{Mod}
\dmo{\PMod}{PMod}
\dmo{\LMod}{LMod}
\dmo{\Diff}{Diff}
\dmo{\Homeo}{Homeo}
\dmo{\dist}{dist}
\dmo\BDiff{BDiff}
\dmo\SO{SO}
\dmo\Sp{Sp}
\dmo\Hom{Hom}
\dmo\SL{SL}
\dmo\rank{rank}
\dmo\sig{sig}
\dmo\Out{Out}
\dmo\Aut{Aut}
\dmo\Inn{Inn}
\dmo\GL{GL}
\dmo\PGL{PGL}
\dmo\Gr{Gr}
\dmo\PSL{PSL}
\dmo\BHomeo{BHomeo}
\dmo\EHomeo{EHomeo}
\dmo\EDiff{EDiff}
\dmo{\PDiff}{PDiff}
\dmo\Disc{Disc}
\dmo\Aff{Aff}
\dmo{\Hol}{Hol}
\dmo{\Map}{Map}

\renewcommand{\hat}{\widehat}
\renewcommand{\bar}{\overline}

\dmo\Teich{Teich}
\dmo\Fix{Fix}
\nc{\pair}[1]{\ensuremath{\left\langle #1 \right\rangle}}
\nc{\abs}[1]{\ensuremath{\left| #1 \right|}}
\nc{\action}{\circlearrowright}
\nc{\abcd}[4]{\ensuremath{\left(\begin{array}{cc} #1 & #2 \\ #3 & #4 \end{array}\right)}}
\nc{\into}\hookrightarrow
\nc{\semidirect}{\rtimes}
\dmo{\Isom}{Isom}
\nc{\normal}{\vartriangleleft}
\dmo{\Vol}{Vol}
\dmo{\im}{Im}
\dmo{\Push}{Push}
\dmo{\Conf}{Conf}
\dmo{\PConf}{PConf}
\dmo{\PB}{PB}
\dmo{\id}{id}
\dmo{\Jac}{Jac}
\dmo{\Pic}{Pic}
\dmo{\ch}{ch}
\dmo{\td}{td}
\dmo{\Stab}{Stab}
\dmo{\Arf}{Arf}
\dmo{\End}{End}
\dmo{\Ext}{Ext}
\dmo{\Gal}{Gal}
\dmo{\lcm}{lcm}
\dmo{\ab}{ab}
\dmo{\opp}{op}
\dmo{\SU}{SU}
\dmo{\OT}{\Omega \mathcal{T}}
\dmo{\OM}{\Omega \mathcal{M}}
\dmo{\PH}{\mathbb{P}\mathcal{H}}
\dmo{\spin}{spin}
\dmo{\even}{even}
\dmo{\odd}{odd}
\dmo{\comp}{\mathcal{H}}
\dmo{\Mgk}{\mathcal{M}_{g, \underline{\kappa}}}
\dmo{\orb}{orb}
\dmo{\AJ}{AJ}
\dmo{\Ck}{\mathsf{C}(\underline{\kappa})}
\dmo{\Int}{Int}
\dmo{\pr}{pr}
\dmo{\lab}{lab}
\dmo{\Sym}{Sym}
\dmo{\Ann}{Ann}
\dmo{\Rad}{Rad}
\dmo{\Ind}{Ind}
\dmo{\Div}{Div}
\dmo{\Res}{Res}
\dmo{\Hur}{Hur}
\dmo{\vcd}{vcd}
\dmo{\codim}{codim}
\dmo{\Emb}{Emb}

\nc{\Span}[1]{\operatorname{Span}(#1)}
\newcommand{\onto}{\twoheadrightarrow}

\renewcommand{\epsilon}{\varepsilon}
\renewcommand{\tilde}{\widetilde}
\renewcommand{\le}{\leqslant}
\renewcommand{\ge}{\geqslant}
\nc{\coloneq}{\mathrel{\mathop:}\mkern-1.2mu=}
\nc{\margin}[1]{\marginpar{\scriptsize #1}}
\nc{\para}[1]{\medskip\noindent\textbf{#1.}}
\definecolor{myblue}{RGB}{102,153, 255}
\definecolor{myred}{RGB}{204,0,0}
\definecolor{mygreen}{RGB}{0,204,0}
\definecolor{myorange}{RGB}{255,102,0}
\definecolor{mypurple}{RGB}{138,43,226}
\nc{\red}[1]{\textcolor{myred}{#1}}
\nc{\blue}[1]{\textcolor{myblue}{#1}}

\dmo{\Heis}{Heis}
\dmo{\FAut}{FAut}

\author{Nick Salter}
\address{Department of Mathematics, University of Notre Dame, 255 Hurley Building, Notre Dame, IN 46556}
\email{nsalter@nd.edu}

\title[Topology of linear systems of hyperelliptic curves in Hirzebruch surfaces]{On the topology of some linear systems of hyperelliptic curves in Hirzebruch surfaces}
\date{November 10, 2026}

\begin{document}

\begin{abstract}
    Following a program suggested by Dolgachev and Libgober in 1981, we study the topology of the locus of smooth sections of linear systems of hyperelliptic curves on Hirzebruch surfaces. Our primary interest is the computation of various flavors of ``enriched monodromy representations''. Within this class of examples, we find that a variety of behaviors occur - in some situations the enriched monodromy faithfully captures the fundamental group, and in others it does not. Our arguments mainly leverage algebro-geometric techniques to compute cohomological obstructions to families of embeddings.
\end{abstract}

\maketitle

\section{Introduction}

We work throughout over $\C$. Let $X$ be a smooth projective variety and $L$ a very ample line bundle on $X$. A program initiated by Dolgachev-Libgober \cite[Section 5]{DL} seeks to understand the topology of the locus of smooth sections in the linear system $\abs{L}$, i.e. the complement of the discriminant, by means of various flavors of {\em topological monodromy representations}. Let $U_L \subset \abs{L}$ denote the discriminant complement, and let $C$ be a smooth section, viewed simultaneously as a basepoint in $U_L$ and an embedded submanifold of $X$. As recalled in \Cref{section:LES}, there are at least three flavors of monodromy to consider. In increasing order of sensitivity, these are the {\em intrinsic monodromy} valued in the mapping class group $\Mod(C) \coloneq \pi_0(\Diff(C))$ of $C$, the {\em relative monodromy} valued in the relative mapping class group $\pi_0(\Diff(X,C))$ of the pair $(X,C)$, and the {\em embedded monodromy} valued in the fundamental group of the ``nonlinear Grassmannian'' $\Gr(C,X)$ of smooth submanifolds of $X$ diffeomorphic to $C$. This is a higher-dimensional analogue of how the braid group (embedded monodromy) gives an enrichment of the symmetric group (intrinsic monodromy).

In this paper, we carry out an analysis of these monodromy representations for certain linear systems on Hirzebruch surfaces. To our knowledge, these are the first complete computations of this kind. 
For $n \ge 0$, let $\bF_n = \bP(\cO_{\P^1} \oplus \cO_{\P^1}(-n))$ denote the $n^{th}$ Hirzebruch surface, equipped with distinguished section $S$ and fiber $F$ (these notions are recalled in \Cref{SS:recap}). Note in particular $\bF_0 \cong \P^1 \times \P^1$; for general readability in the introduction we will use the latter notation, but will use the former in the body of the paper. The smooth sections of the  linear system $\abs{2S + (n+g+1)F}$ are hyperelliptic curves of genus $g$. We write $U_{n,g}$ for the associated locus of smooth sections; it is nonempty if and only if $n \le g+1$.   

We focus our attention on two classes of examples at opposite ends of this spectrum, the cases $n = 0$ and $n = g$. In each of these cases, we find (cf. \Cref{main:bundle}) that the coarsest (intrinsic) monodromy is non-injective, warranting a study of the finer monodromy maps suggested by Dolgachev-Libgober. Intriguingly, we find that the cases $n = 0$ and $n = g$ exhibit very different behavior.

\begin{maintheorem}\label{main:P1xP1}
    For all $g \ge 2$, the embedded monodromy homomorphism $\epsilon_*: \pi_1(U_{0,g}) \to \pi_1(\Gr(C, \bP^1 \times \bP^1))$ is injective, while the reduction to the relative monodromy $\bar \epsilon_*: \pi_1(U_{0,g}) \to \pi_0(\Diff(\bP^1 \times \bP^1, C))$ has kernel $\Z \times (\Z/2\Z)^2$, and the intrinsic monodromy map $\rho_*: \pi_1(U_{0,g}) \to \Mod(\Sigma_g)$ has kernel $\Heis(H_1(C;\Z)) \times (\Z/2\Z)^2$.
\end{maintheorem}
In the above, $\Heis(H_1(C;\Z))$ denotes the {\em Heisenberg group} associated to $H_1(C;\Z)$, i.e. the cyclic central extension of $H_1(C;\Z)$ with extension class given by the algebraic intersection form. Thus we find that in the case $n = 0$, the entire fundamental group can be detected topologically via a suitably sensitive monodromy representation. This is in contrast with what happens in the case $n = g$: we find that even the finest monodromy representation fails to capture the whole fundamental group.

\begin{maintheorem}\label{main:Fg}
    For all $g \ge 2$, the kernel of the embedded monodromy $\epsilon_*: \pi_1(U_{g,g}) \to \pi_1(\Gr(C, \bF_g))$ contains a central extension of the third term in the lower central series for $\pi_1(C)$. The image is an extension of the hyperelliptic mapping class group by a two-step nilpotent group given as a central extension of $H_1(C;\Z)$ by an abelian group of rank one.
\end{maintheorem}

A more precise formulation of this theorem is given in \Cref{thm:preciseFg}.
An intriguing byproduct of this latter result concerns the fundamental group of the diffeomorphism groups of the smooth manifolds underlying $\bF_g$. Recall that as a smooth manifold, $\bF_g$ is an $S^2$-bundle over $S^2$, diffeomorphic to $S^2 \times S^2$ if $g$ is even, and to $\bP^2 \# \bar{\bP^2}$ if $g$ is odd. The proof of \Cref{main:Fg} implies the following result (the notation $\Diff_0$ denotes the identity component of the diffeomorphism group).

\begin{maincor}\label{cor:diff}
The groups $\pi_1(\Diff_0(S^2 \times S^2))$ and $\pi_1(\Diff_0(\bP^2 \# \bar{\bP^2}))$ each contain an element of infinite order.
\end{maincor}

At least in the case of $S^2 \times S^2$ this is not the first example of its kind (see {\bf Related work} below), but we felt the underlying mechanism was of sufficient interest as to warrant mention.\\

These results are built off of an explicit topological description of the discriminant complements as orbi-bundles, leading in particular to a description of the fundamental groups. The result below is a modest strengthening of work of Bergstrom-Zheng \cite{bergstromzheng}, who constructed moduli stacks of hyperelliptic curves in Hirzebruch surfaces and studied their relationship with the moduli stack $\cH_g$ of hyperelliptic curves, but stopped short of giving a description of the fundamental groups of the fibers.

\begin{maintheorem}\label{main:bundle}
    For all $g \ge 2$, the classifying map $U_{0,g} \to \cH_g$ is an orbifold fiber bundle with fiber $(\Pic^{g+1}(C) \setminus \Theta') \times (\PGL_2(\C))^2$, where $\Theta \subset \Pic^{g-1}(C)$ is a theta-divisor, $H \in \Pic^2(C)$ is the hyperelliptic line bundle, and $\Theta' = H \otimes \Theta$. Consequently there is a short exact sequence
    \[
    1 \to \Heis(H_1(C;\Z))\times (\Z/2\Z)^2 \to \pi_1(U_{0,g}) \to \Delta_g \to 1,
    \]
    where $\Delta_g \coloneq \pi_1^{orb}(\cH_g) \le \Mod(\Sigma_g)$ denotes the hyperelliptic mapping class group. 
    
    Similarly, the classifying map $U_{g,g} \to \cH_g$ is an orbifold fiber bundle with fiber $E_C$ a principal $\Aut(\bF_g)$-bundle over $C$, leading to the short exact sequence
    \[
    1 \to \pi_1(E_C) \to \pi_1(U_{g,g}) \to \Delta_g \to 1.
    \]
    Moreover $\pi_1(E_C)$ is realized as a nonsplit central extension
    \begin{equation}\label{eq:pi1EC}
        1 \to \pi_1(\Aut(\bF_g)) \to \pi_1(E_C) \to \pi_1(C) \to 1,
    \end{equation}
    with $\pi_1(\Aut(\bF_g))$ realized explicitly as the abelian group $\Z\pair{v,w}/(v^g = w^2)$; the Euler class of \eqref{eq:pi1EC} is $v^2 \in H^2(C; \pi_1(\Aut(\bF_g)))$.
\end{maintheorem}

\para{Ideas of proof} As mentioned above, the proof of \Cref{main:bundle} is largely contained in work of Bergstrom-Zheng, though we give a self-contained and somewhat different proof that is more suited to our purposes. The key idea is that an embedding $C \into \bF_n$ of a hyperelliptic curve into a Hirzebruch surface is governed by a line bundle on $C$ with certain cohomological properties. The embedding arises by pushing this line bundle forward along the hyperelliptic map to obtain a rank-2 vector bundle on $\P^1$, leading to a fiberwise embedding of $C$ into the associated projective bundle. For fixed $C$ and $n$, the locus of suitable line bundles can be further recast as a space of divisors, leading to an understanding of the fundamental group. This is carried out in \Cref{section:hypintohirz}.

To prove that the embedded monodromy (arising as $\pi_1$ of the natural map $\epsilon: U_{0,g} \to \Gr(C, \bF_0)$) detects $\pi_1(U_{0,g})$, we must construct a certificate that a given loop of embedded surfaces (arising from a loop in $U_{0,g}$) is not null-homotopic as a loop of smooth submanifolds, i.e. as a loop in $\Gr(C,\bF_0)$. Such a null-homotopy would give rise to an extension of the fiberwise embedding $C \times S^1 \to \bF_0 \times S^1$ over a disk bounding $S^1$. We find that the characteristic classes of the normal bundle of such a hypothetical embedding obstruct its existence - this argument is realized in \Cref{lemma:thetameridianinforder}.

As discussed in \Cref{section:Fg}, the mechanism of proof in \Cref{main:Fg} is to show that the map $\epsilon: E_C \to U_{g,g} \to \Gr(C,\bF_g)$ factors through an enlargement $\cE_C$ of $E_C$ whose fundamental group is two-step nilpotent. Such an extension exists for all $n$ (cf. \Cref{section:universalfamily}), but in the case $n = 0$, the relevant group is already two-step nilpotent and so this does not obstruct injectivity. Thus, \Cref{main:bundle} offers a philosophical explanation for the difference between the cases $n = 0$ and $n = g$ studied in \Cref{main:P1xP1,main:Fg}: the kernel term in the case $n = g$ is ``too nonabelian''. 

Finally, \Cref{cor:diff} on the fundamental group of the diffeomorphism groups is a corollary of the proof of \Cref{main:Fg}. As discussed in \Cref{section:LES}, the various monodromy representations we study are related by homotopy long exact sequences for fibrations of diffeomorphism groups, {\em \`a la } the classic work of Cerf and Palais. Via this exact sequence, elements in the kernel of the relative monodromy that are detected by the embedded monodromy then lift to elements in the fundamental group of the ambient diffeomorphism group.

\para{Related work} We have already commented on the work of Bergstrom-Zheng \cite{bergstromzheng} and its role in \Cref{main:bundle}, as well as the pioneering work of Dolgachev-Libgober \cite{DL} which served as the inspiration for the project as a whole. The Dolgachev-Libgober program has been taken up by various authors over the years - we survey some of the notable results in the area. 

In \cite{carlsontoledo}, Carlson-Toledo study a secondary homological monodromy map present for linear systems of projective hypersurfaces, obtained by taking the branched cover of $\P^n$ along the hypersurface as it varies in the linear system. Evidently this factors through the embedded monodromy representation as we have defined it here, though it is vastly more computable. Carlson-Toledo use this to show that the kernel of the ``classical'' monodromy (acting on the homology of the hypersurface) is quite large. In \cite{saltermonker}, the author extends Carlson-Toledo's result to the topological setting, showing that the intrinsic monodromy representation for the family of smooth plane curves is large (in the same sense as Carlson-Toledo) for all degrees $d \ge 4$ (previously, the case $d = 3$ was analyzed by Dolgachev-Libgober, and the case $d=4$ was described completely in the work of R. Harris \cite{harris}). 

Another crucial antecedent is the work of Bolognesi-L\"onne \cite{bolognesilonne}, who carried out a study of linear systems of {\em trigonal} curves on Hirzebruch surfaces. They were able to arrive at a complete description of the intrinsic monodromy, and in the case of the ``maximal Maroni stratum'', obtained an explicit presentation for the fundamental group; they did not study the relative or embedded monodromy maps. 

The intrinsic monodromy of a linear system of curves in a simply-connected surface $X$ (valued in a mapping class group $\Mod(C)$) is by now well-understood for sufficiently ample line bundles. Work \cite{nickishan} of the author with I. Banerjee shows that under a certain ampleness hypothesis, the monodromy group $\rho_*(\pi_1(U_L)) \le \Mod(C)$ of a linear system $\abs{L}$ is a finite-index subgroup known as an {\em $r$-spin mapping class group} - the $r$-spin structure arises as a maximal root of the adjoint line bundle. But the work of \cite{saltermonker} shows that the kernel of $\rho_*$ is always infinite in these cases. It is unknown how much more is detected by relative or embedded monodromy.

Even the case of linear systems of divisors on curves is quite rich. Here, the discriminant complement of a linear system of degree $d$ can be identified with a locus in the configuration space $\Conf_d(C)$ of $d$ points in $C$ (as smoothness is equivalent to the divisor consisting of $d$ distinct points). Formulated in our language, Dolgachev-Libgober sketched a proof that when $d$ is sufficiently large, the embedded monodromy (valued here in the surface braid group $\pi_1(\Conf_d(C))$) is injective. They also identified the image: under the Abel-Jacobi map, the locus in $\Conf_d(C)$ coming from smooth sections of a fixed line bundle $L$ is {\em constant}, and so the fundamental group lies in the ``simple braid group'', the kernel of the Abel-Jacobi map $\lambda_*: \pi_1(\Conf_d(C)) \to H_1(C;\Z)$. Dolgachev-Libgober's sketch-proof was refined and confirmed in work of Shimada \cite{shimada}.

Beyond Carlson-Toledo's work, another notable contribution to the story in higher dimensions is the book chapter of Hain \cite{hain}, which lays out a number of interesting results, techniques, and questions. 
And very recently, Aumonier has made a number of very interesting contributions to the study of the topology of discriminant complements, in particular applying more sophisticated methods from algebraic topology ($h$-principles, scanning maps, Goodwille calculus, etc.) - see \cite{aumonier1,aumonier4,aumonier2}.

With regards to the literature on the homotopy type of diffeomorphism groups of $4$-manifolds, the landmark result here is Watanabe's disproof of the Smale conjecture \cite{watanabe}, which in particular produces infinite-order elements in $\pi_1(\Diff(S^4))$. In personal communication, D. Ruberman has explained to us that J. Lin's work \cite{lin} also implicitly produces infinite-order elements in $\pi_1(\Diff(S^2 \times S^2))$, which can be detected using a similar Pontryagin class obstruction. 

\para{Speculations and discussion} The non-injectivity result in \Cref{main:Fg} comes as something of a disappointment to the author, who had long been hoping that the embedded monodromy would always be sensitive enough to faithfully capture the fundamental group of the discriminant complement of an arbitrary linear system, at least for curves in surfaces. The results here do shed some light on the mechanisms underlying (non)injectivity, and the structure of the fundamental groups of discriminant complements more generally, which are worth commenting briefly on. 

Let $X$ be a smooth projective surface and $L$ a very ample line bundle on $X$ with discriminant complement $U_L$. There is always a classifying map $\rho: U_L \to \cM_g$, and the induced map on orbifold fundamental group is the intrinsic monodromy. There are then two separate mechanisms that contribute to $\ker(\rho_*)$. The first is that the image $\rho(U_L) \subset \cM_g$ can be embedded in a non $\pi_1$-injective way. This is the case, for example, for the family of smooth plane quartics: the image in $\cM_g$ is the complement of the hyperelliptic locus, and the fundamental group of this was shown by R. Harris \cite{harris} to be an extension of $\Mod(\Sigma_3)$ by a free group of countably infinite rank. 

The structural result of \Cref{main:bundle} highlights a separate mechanism giving rise to elements of $\ker(\rho_*)$: loops in the fiber of $\rho$. This is present in the classical examples studied in Dolgachev-Libgober's original work, where the fundamental group of certain linear systems of elliptic curves are realized as extensions of $\SL_2(\Z)$ by finite groups which arise as isotropy groups for the action of the automorphism group of the background space. The linear systems of hyperelliptic curves in Hirzebruch surfaces studied in this paper are interesting from this perspective precisely because a given $C$ has a large and topologically rich space of embeddings into a fixed $\bF_n$.

This perspective also clarifies one of the fundamental challenges in understanding a given $\pi_1(U_L)$ (focusing attention here on linear systems of curves in surfaces). Namely, even if it were true in some greater generality that $\rho: U_L \to \cM_g$ is some kind of orbifold fibration over its image, this still leaves the problem of determining the fundamental group of the base space $\rho(U_L)$. The examples in this paper circumvent this issue because the hyperelliptic locus in $\cM_g$ is well-understood, but as mentioned above, even in the relatively simple setting of plane quartics, this becomes a very delicate problem. 

The question of the structure of the filtration of $\ker(\rho_*)$ given by $\ker(\epsilon_*) \le \ker(\bar \epsilon_*) \le \ker(\rho_*)$ is an intriguing one. The main mechanism driving \Cref{lemma:thetameridianinforder} (which shows that the containment $\ker(\epsilon_*) < \ker(\bar \epsilon_*)$ is strict in the case of $\bF_0$) seems to be a special feature of linear systems on Hirzebruch surfaces: a family of smooth embeddings can degenerate to a {\em smooth} embedding into a surface that is diffeomorphic but not biholomorphic to the ambient space at the generic point. It is not clear to the author whether more examples of this type exist.

\para{Organization of paper} In \Cref{section:LES}, we review the construction of the monodromy maps. In \Cref{section:hypintohirz}, we review the theory of embeddings of hyperelliptic curves in Hirzebruch surfaces and give a proof of \Cref{main:bundle}. In \Cref{section:universalfamily}, we construct a ``smooth universal family'' of embeddings of a fixed hyperelliptic $C$ into Hirzebruch surfaces; this construction is essential for our study of embedded monodromy. In \Cref{section:F0}, we study the case of hyperelliptic curves in $\bF_0$, proving \Cref{main:P1xP1}, and in \Cref{section:Fg}, we prove \Cref{main:Fg}.

\para{Acknowledgments} I thank Ishan Banerjee, Sander Kupers, Danny Ruberman, and Tadayuki Watanabe for helpful bibliographic input. This work was supported in part by NSF CAREER grant no. DMS-2338485.

\para{AI usage} This project was carried out with extensive assistance from ChatGPT - the bulk of the technical ideas of the paper originated on the machine side of the equation. In an attempt to minimize the deficiencies evident in AI-driven work, I have written this paper the old-fashioned way (completely free of AI assistance), after a period of studying the arguments supplied by the machine, reformulating them to suit my expositional preferences and to better reflect my pre-existing background knowledge. The final draft of the paper was run through the machine to check for errors, but again all edits were made by hand. Every effort has been made to track down appropriate references and to identify precedents for the ideas used in the paper, but I fully concede the possibility that this may be incomplete; I warmly welcome bibliographic suggestions.

\section{Fiber sequences for diffeomorphism groups and embedding spaces}\label{section:LES}

Here we recollect some basic results on structural relationships between diffeomorphism groups and embedding spaces, in order to give a precise formulation of the various monodromy maps under consideration. For smooth manifolds $X,Y$, we let $\Emb(Y,X)$ denote the space of smooth embeddings $f: Y \to X$, endowed with the $C^\infty$ topology. When $Y = X$, we identify $\Emb(Y,X)$ with the diffeomorphism group $\Diff(X)$. Finally, for a closed submanifold $C \subset X$, we let $\Gr(C,X)$ denote the quotient $\Emb(C,X)/\Diff(C)$, the ``nonlinear Grassmannian'' of smoothly embedded submanifolds of $X$ that are diffeomorphic to $C$. 

The following result lays the technical foundations for the construction of our monodromy maps. A proof can be found in \cite[Lemma 2.1a,b]{hallervizman}; see also the treatment in \cite[Theorem 44.1]{michor} (the original work of Palais \cite{palais} and Cerf \cite[Appendice]{cerf} both treat only the ``parametrized'' setting, replacing $\Gr(C,X)$ with the space of embeddings $\Emb(C,X)$ and so cannot be used verbatim). 

\begin{theorem}\label{theorem:difffibr}
    Let $X$ be a smooth compact manifold without boundary, and let $C \subset X$ be a closed submanifold, also without boundary. Then the orbit map
    \[
    \Diff(X) \to \Gr(C,X)
    \]
    is a locally trivial fibration over its image. The fiber over the distinguished point of $\Gr(C,X)$ represented by $C \subset X$ is the group $\Diff(X,C)$ of diffeomorphisms of the pair $(X,C)$. 
    In addition, the quotient $\Emb(C,X) \to \Gr(C,X)$ is also a locally trivial fibration, with fiber given by $\Diff(C)$.
\end{theorem}

Using this result, we can construct the various monodromy maps considered in the paper. Let $X$ be a smooth algebraic surface, $L$ a line bundle on $X$, and $U_L$ the complement of the discriminant locus in the linear system $\abs{L}$, equipped with basepoint $C$ (viewed simultaneously as a smooth section of $L$ and an embedded submanifold of $X$). Then there is a natural continuous map
\[
\epsilon: U_L \to \Gr(C,X)
\]
and hence an induced map $\epsilon_*: \pi_1(U_L) \to \pi_1(\Gr(C,X))$. We call this the {\em embedded monodromy}. Composing this with the boundary map $\delta: \pi_1(\Gr(C,X)) \to \pi_0(\Diff(X,C))$ in the long exact sequence for the fibration $\Diff(X) \to \Gr(C,X)$, we obtain the second monodromy map, called the {\em relative monodromy}
\[
\bar \epsilon_*: \pi_1(U_L) \to \pi_0(\Diff(X,C)).
\]
The {\em intrinsic monodromy} $\rho_*: \pi_1(U_L) \to \Mod(C) = \pi_0(\Diff^+(C))$ can then be obtained by composing $\bar \epsilon_*$ with $\pi_0$ of the natural restriction $\Diff(X,C) \to \Diff(C)$ (or equivalently, by taking the induced map on orbifold fundamental group of the classifying map $\rho: U_L \to \cM_g$).

In the sequel it will be convenient to replace the target $\Gr(C,X)$ of the space-level map $\epsilon$ defining the embedded monodromy with the richer space $\Emb(C,X)$ (i.e. replace embedded submanifolds $C \subset X$ with embeddings $i: C \to X$). Below we identify when this is possible.

\begin{lemma}\label{lemma:kerrholifts}
    The kernel $\ker(\rho_*) \le \pi_1(U_L)$ of the intrinsic monodromy admits a lift $\tilde \epsilon_*: \ker(\rho_*) \to \pi_1(\Emb(C,X))$.
\end{lemma}
\begin{proof}
The long exact sequence for the fibration $\Emb(C,X) \to \Gr(C,X)$ includes the segment 
\[
\pi_1(\Diff_0(C)) \to \pi_1(\Emb(C,X)) \to \pi_1(\Gr(C,X)) \to \pi_0(\Diff(C)).
\]
As we assume $g \ge 2$, the Earle-Eells theorem \cite{earleeells} states that $\pi_1(\Diff_0(C)) = 0$. Exactness now identifies $\pi_1(\Emb(C,X))$ with the kernel of the map $\pi_1(\Gr(C,X)) \to \pi_0(\Diff(C))$. As $\rho_*$ factors through this map, the result follows. 
\end{proof}

\section{Linear systems on hyperelliptic curves and embeddings into Hirzebruch surfaces}\label{section:hypintohirz}
    \subsection{Hyperelliptic curves and Hirzebruch surfaces}\label{SS:recap}
    Here we fix notation and recall the rudiments of the theories of hyperelliptic curves and Hirzebruch surfaces.
    
    \para{Hyperelliptic curves} Throughout, $g \ge 2$ will be fixed, and $C$ will denote a smooth hyperelliptic curve of genus $g$. A {\em hyperelliptic map} $h: C \to \P^1$ is any meromorphic function of degree $2$; such $h$ is unique up to the action of $\Aut(\P^1) = \PGL_2(\C)$. The {\em hyperelliptic line bundle} $H \in \Pic^2(C)$ is defined via 
    \[
    H:= h^*(\cO_{\P^1}(1)),
    \]
    and the hyperelliptic involution is denoted 
    \[
    \iota: C \to C.
    \]

    \para{Hirzebruch surfaces} We will assume familiarity with the basic theory of Hirzebruch surfaces at the level of \cite[pp. 514-520]{GH}; let us fix our notational conventions. For $n \ge 0$, the Hirzebruch surface $\bF_n$ is defined as the projective bundle\footnote{Here and throughout, we use the quotient convention for projective spaces: $\bP(V)$ parametrizes rank-one {\em quotients} of $V$.} associated to $\cO \oplus \cO(-n)$ on $\P^1$. Recall the basic fact that projectivization ignores twisting: $\bP(\cO(a) \oplus \cO(a-n)) \cong \bP((\cO\oplus \cO(-n))\otimes \cO(a)) \cong \bP(\cO \oplus \cO(-n))$. In particular, the isomorphism type of an abstract projective bundle $\bP(V)$ (for $V$ a rank-$2$ vector bundle on $\CP^1$) is determined solely by the difference $p-q$ under an isomorphism $V \cong \cO(p) \oplus \cO(q)$.
    
    Projection gives a {\em ruling map} 
    \[
    p_n: \bF_n \to \P^1,
    \]
    with fiber class $F$ an irreducible rational curve satisfying $F^2 = 0$. The section $s: \P^1 \to \bF_n$ given by the distinguished quotient $\cO \oplus \cO(-n) \to \cO(-n)$ gives a second irreducible rational curve called the {\em distinguished section}. For $n > 0$ it is the unique irreducible rational curve satisfying $S^2 = -n$ (this can be computed by observing that its normal bundle is isomorphic to $\cO(-n)$). The description of the intersection form on $\bF_n$ is then completed by observing that $SF = 1$. We have explicitly allowed $n = 0$ throughout; in this case, $\bF_n$ is trivial as a $\P^1$-bundle, so that $\bF_0 \cong \P^1 \times \P^1$. The existence of a second ruling here will play a relatively minor role in our story.

    Concerning the {\em topology} of Hirzebruch surfaces, we recall the fact that there are precisely two isomorphism classes of oriented $S^2$-bundles over $S^2$, and that $\bF_n$ is diffeomorphic to $\bF_m$ if and only if $n \equiv m \pmod 2$ (see \cite[Theorem 3.4.8]{gompfstipsicz}).

    \subsection{Families of embeddings into Hirzebruch surfaces}
    Here we consider the linear systems $\abs{2S + (n+ g + 1)F}$ on $\bF_n$, for $n \ge 0$ and $g \ge 2$ arbitrary. For $C \in \abs{2S + (n+g+1)F}$ a smooth section, restricting the projection along the ruling yields a map $h: C \to \P^1$ of degree $2$, showing that $C$ is hyperelliptic; adjunction shows that $C$ has genus $g$. In this section, we review the work of Bergstrom-Zheng \cite{bergstromzheng}, which constructs the moduli space of hyperelliptic curves embedded in $\bF_n$ in a particular linear system and implicitly gives the orbi-bundle description of the linear system underlying our \Cref{main:bundle}. Bergstrom-Zheng obtain their results via an embedding of $\bF_n$ into a weighted projective space. In \Cref{lemma:hyptohirz} and \Cref{prop:Ubundle}, we present here a different and somewhat more geometric perspective on their results (also known to Bergstrom-Zheng - see \cite[Remark 5.2]{bergstromzheng}). This is a key ingredient in the proof of \Cref{main:bundle}. It remains to give a description of the fundamental group of the fiber, which was determined by Nori \cite[Corollary 2.9]{nori} in the case $n = 0$ (cf. \Cref{lemma:thetaheis}), and which is worked out in the case $n = g$ below in \Cref{SS:ECng}.
    
    \begin{lemma}[Cf. \cite{bergstromzheng}, Proposition 5.3]\label{lemma:hyptohirz}
        Let $1 \le d \le g+1$, and let $L \in \Pic^d(C)$ be given. Define $a\in \Z$ as the least integer such that $h^0(L \otimes H^a) > 0$, and set $n = g+1-d-2a$. Then $L$ determines an embedding $i: C \into \bF_{n}$, realizing $i(C)$ as a smooth section in the linear system $\abs{2S + (n+ g+1)F}$. 

        Conversely, for $0 \le n \le g$, up to the action of $\Aut(\bF_n)$, a smooth section of $\abs{2S+(n+ g+1)F}$ is uniquely determined by a smooth curve $C \in \cH_g$ together with a line bundle $L \in \Pic^{g+1-n}(C)$ such that $h^0(L) > 0$ but $h^0(L\otimes H^{-1}) = 0$.
    \end{lemma}
    \begin{proof}
        That $L$ determines a morphism of $C$ into a ruled surface is a standard construction which we only sketch here; see \cite[Proposition 7.12]{hartshorne}. 
        Pushing $L$ forward along the finite morphism $h: C \to \bP^1$ determines a vector bundle $h_*L$ of rank $\deg(h) = 2$ on $\bP^1$. The evaluation quotient $h^*h_*L \onto L$ gives a morphism
        \[
        i: C \to \bP(h_* L)
        \]
        compatible with the projections $h: C \to \bP^1$ and $\pi: \bP(h_* L) \to \bP^1$, and it is straightforward to see that this separates points and tangent vectors and so determines an embedding.
        
        It remains to identify $\bP(h_* L)$ with a specific Hirzebruch surface.
        By Riemann-Roch for vector bundles,
        \[
        \chi(\bP^1,h_*L) = \deg(h_*L) + 2,
        \]
        and as $h$ is finite and hence affine, $h^i(h_*L) = h^i(L)$ for all $i$. Thus by ordinary Riemann-Roch,
        \[
        1-g + d = \chi(C,L) = \chi(\bP^1,h_*L) = \deg(h_*L) + 2,
        \]
        showing 
        \[
        \deg(h_* L) =  d - 1 - g.
        \]
        By the splitting theorem for vector bundles on $\bP^1$, it follows that
        \[
        h_*L \cong \cO(p) \oplus \cO(q)
        \]
        for integers $p \ge q$ satisfying $p+q = d - 1 - g$; consequently we find 
        \[
        \bP(h_* L) \cong \bF_{2p+g+1-d}.
        \]

        It remains to determine $p$. To do so, recall that $H = h^*\cO_{\P^1}(1)$, so that for any integer $r$,
        \[
        h_*(L\otimes H^r) \cong \cO(p+r) \oplus \cO(q + r).
        \]
        As $p \ge q$ by assumption, and since 
        \[
        h^0(L \otimes H^r) = h^0(h_*(L \otimes H^r)) = h^0(\cO(p+r) \oplus \cO(q+r)),
        \]
        we see that $h^0(L \otimes H^r) > 0$ if and only if $p+r \ge 0$. This shows that $p = -a$, where $a$ is defined as above as the minimal integer such that $h^0(L \otimes H^a) > 0$. 

        To determine the divisor class of $i(C)$, it is first clear that $i(C) \cdot F = 2$, so that $[i(C)] = 2S + mF$ for some $m$. We claim that $i(C) \cdot S = g+1-n = d+2a$, and moreover that $i(C) \cap S$ lies in the linear system $\abs{L \otimes H^a}$. 
        To see this, recall that the embedding $i: C \to \bP(h_*L)$ is given by the evaluation quotient $h^*h_*L \onto L$. After identifying $h_*L \cong \cO(-a) \oplus \cO(-a-n)$, restrict this quotient to the sub-bundle $h^*(\cO(-a))$ to obtain a morphism of line bundles
        \[
        h^*(\cO(-a)) \to L,
        \]
        yielding a section $s$ of $\Hom(h^*(\cO(-a)),L) \cong L \otimes h^*(\cO(a))$, say with divisor $D_s$. As $h^*(\cO(1)) = H$ by definition, this shows that $D_s \in \abs{L \otimes H^a}$. On the other hand, after identifying $\bP(h_*L) \cong \bF_n$, the negative section $S$ on $\bF_n$ is given by the quotient $\cO \oplus \cO(-n) \onto \cO(-n)$, and so before making this identification, $S$ is given by the quotient annihilating $\cO(-a)$. Thus $i(C) \cap S$ is precisely the vanishing locus of $s$, i.e. $D_s$.

        To prove the converse, suppose that $0 \le n \le g$, and let $C \in \abs{2S + (n+ g+1)F}$ be a smooth section. Then $C \cap S$ is an effective divisor on $C$ of degree $C \cdot S = g+1-n$. Moreover, $C \cap S$ contains no conjugate pair $p + \iota(p)$: any such pair arises as the intersection $C \cap F$ with some fiber of the ruling, and so if $C \cap S$ contained such a pair, it would violate the intersection calculation $S \cdot F = 1$. Taking $L = \cO_C(C \cap S)$, this establishes the converse.
    \end{proof}
    
The previous result shows that embeddings of a hyperelliptic curve $C$ into a ruled surface are parametrized by suitable line bundles $L$ on $C$. To translate this into an understanding of the smooth sections in $\abs{2S + (n+g+1)F}$, we need to specify not just $L \in \Pic^d(C)$, but also an identification $\P(h_*(L)) \cong \bF_n$. Accordingly, with $n,g$ implicitly fixed, define
    \[
   F_C \coloneq \{L \in \Pic^{g+1-n}(C) \mid h^0(L) > 0,\ h^0(L \otimes H^{-1}) = 0\},
   \]
    and let 
    \begin{equation}
        \label{eq:EC}
        E_C = \{(L,h, \Phi) \mid L \in F_C,\ h: C \to \P^1\mbox{ of degree $2$, } \Phi: \bP(h_*L) \cong \bF_n\mbox{ a bundle iso.}\}.
    \end{equation}
    With an eye towards giving a more concrete description of $F_C$, we define
    \[
    \Sym^{g+1-n}(C)^\circ \subset \Sym^{g+1-n}(C)
    \]
    as the locus of effective divisors on $C$ not containing any hyperelliptic conjugate pair $p + \iota(p)$. With these notations established, we can continue our recap of Bergstrom-Zheng's work. In order to treat the case $n = 0$ uniformly, define 
    \[
    \Aut'(\bF_n) \le \Aut(\bF_n)
    \]
    as the subgroup that preserves each ruling - for $n \ge 1$, $\Aut'(\bF_n) = \Aut(\bF_n)$, while for $n = 0$ it is a subgroup of index two.

    \begin{proposition}\label{prop:Ubundle}
        For all $g \ge 2$ and $0 \le n \le g$, the classifying map
        \[
        \rho: U_{n,g} \to \cH_g
        \]
        is an orbifold fiber bundle with fiber isomorphic to $E_C$. For $1 \le n \le g$, the fiber $E_C$ is isomorphic to a principal $\Aut'(\bF_n)$-bundle over $\Sym^{g+1-n}(C)^\circ$, while for $n = 0$, it is isomorphic to $\Sym^g(C \setminus \{w\})^\circ \times \Aut'(\bF_0)$, where $w \in C$ is a Weierstrass point.
    \end{proposition}
    
\begin{proof}
   To show that $\rho: U_{n,g} \to \cH_g$ is an orbifold fiber bundle, it suffices to show that its pullback to a manifold cover of $\cH_g$ is a genuine fiber bundle. We will do so by tacitly pulling back to (a connected component of) the hyperelliptic locus $\tilde \cH_g$ in Teichm\"uller space, i.e. equipping $C \in \cH_g$ with a marking.
   
   By \Cref{lemma:hyptohirz}, the fiber (in the orbifold sense) of $\rho: U_{n,g} \to \cH_g$ over $C \in \cH_g$ is a principal $\Aut'(\bF_n)$-bundle over the locus 
   \[
   F_C = \{L \in \Pic^{g+1-n}(C) \mid h^0(L) > 0,\ h^0(L \otimes H^{-1}) = 0\}.
   \]
   It is easy to see that $E_C$ is likewise a principal $\Aut'(\bF_n)$-bundle over $F_C$, and that there is an isomorphism of bundles $E_C \cong \rho^{-1}(C)$.
   
   We next identify $F_C$ with an explicit space of divisors. For $1 \le n \le g$, geometric Riemann-Roch (or the analysis in the proof of \Cref{lemma:hyptohirz}) shows that $h^0(L) = 1$, so that $F_C$ can be identified with the locus of effective divisors
   \[
   \Sym^{g+1-n}(C)^\circ \coloneq \{D \in \Sym^{g+1-n}(C) \mid D \mbox{ contains no conjugate pair $p + \iota(p)$}\}.
   \]
   
   For $n = 0$, there is also an identification of $F_C$ with a space of effective divisors. Here, $h^0(L) = 2$ and it is easy to see that $L$ is basepoint free, so that $h^0(L \otimes \cO_C(-w)) = 1$ for any $w \in C$. Taking $w$ a Weierstrass point, we therefore identify $F_C \subset \Pic^{g+1}(C)$ with
   \[
   \Sym^{g}(C \setminus \{w\})^\circ \coloneq \{D \in \Sym^{g}(C \setminus\{w\}) \mid D \mbox{ contains no conjugate pair $p+\iota(p)$}\}.
   \]
   
   This concrete description of $F_C$ shows that the family of spaces $F_C$ forms a locally trivial fiber bundle over $\tilde \cH_g$. Indeed, a local $\iota$-equivariant trivialization of $\tilde \cH_g$ induces a trivialization of the bundles of symmetric powers; $\iota$-equivariance implies that these restrict to a trivialization of $\Sym^{g+1-n}(C)^\circ$ and $\Sym^{g}(C \setminus \{w\})^\circ$.

   It remains to show that in the case $n = 0$, the principal bundle $E_C/F_C$ is a product. To do so, it suffices to exhibit a section $s: F_C \to E_C$, i.e. to describe a canonical embedding $i: C \to \bF_0$ given only the data of an effective divisor $D \in F_C$. As $\Aut'(\bF_0) \cong (\PGL_2(\C))^2$ acts by automorphisms of each of the rulings, this is equivalent to normalizing the projections $C \to \bP^1$ onto each of the factors. In one direction, this projection coincides with the hyperelliptic map $h: C \to \bP^1$. This can be normalized by choosing three Weierstrass points $w_1, w_2, w_3$ on $C$ and insisting that $h$ send them to $0,1, \infty$, respectively.

    It is necessary to also normalize the projection $\pi': C \subset \bP^1 \times \bP^1 \to \bP^1$ onto the second factor. We do so by normalizing the $2$-jet of $\pi'$ at the Weierstrass point $w_1$, as follows. There are precisely two choices of local coordinate $z$ on $C$ near $w_1$ in which $h: C \to \bP^1$ is given as $z^2$; such a choice can be made locally consistently as $C$ varies. There is then a unique coordinate on $\bP^1$ in which $\pi'(z) = z + O(z^3)$. Together with the normalization of $h$ previously discussed, this gives the required section of $E_C/F_C$.     
\end{proof}

In the case $n = 0$, the fundamental group of $F_C$ can be determined using work of Nori \cite{nori}.

\begin{lemma}\label{lemma:thetaheis}
    Let $C$ be a hyperelliptic curve, and let $w \in C$ be a Weierstrass point. Then the space $F_C=\Sym^g(C\setminus\{w\})^\circ$ is isomorphic to $\Pic^{g-1}(C) \setminus \Theta$, the complement of a theta divisor, and $\pi_1(\Pic^{g-1}(C) \setminus \Theta) \cong \Heis(H_1(C;\Z))$.
\end{lemma}
\begin{proof}
    The proof of \Cref{prop:Ubundle} shows that $\Sym^g(C\setminus\{w\})^\circ$ is isomorphic to the locus
    \[
    \{L \mid h^0(L) > 0,\ h^0(L \otimes H^{-1}) = 0\} \subset \Pic^{g+1}(C).
    \]
    The map $H^{-1} \otimes: \Pic^{g+1}(C) \to \Pic^{g-1}(C)$ identifies this with $\Pic^{g-1}(C) \setminus W_{g-1}$, which by Riemann's theorem can be identified with the complement of $\Theta$.

    By \cite[IV, Corollary 4.5]{ACGH}, the singular locus of $W_{g-1}(C)$ is $W_{g-1}^1(C)$, and by the geometric Riemann-Roch theorem (cf. \cite[p. 13]{ACGH}), $\dim(W_{g-1}^1(C)) = g-3$. According to \cite[Corollary 2.9]{nori}, if the singular support of $\Theta$ is codimension at least three in $\Pic^{g-1}(C)$, then $\pi_1(\Pic^{g-1}(C) \setminus \Theta) \cong \Heis(H_1(C;\Z))$; the claim follows.  
\end{proof}

\begin{proof}[Proof of \Cref{main:bundle}, case $n = 0$]
    $\cH_g$ is an aspherical orbifold for all $g \ge 2$ \cite[Proposition 1]{gonzalezdiezharvey}, so that the long exact sequence in homotopy for $\rho: U_{0,g} \to \cH_g$ given by \Cref{prop:Ubundle} degenerates to the short exact sequence of orbifold fundamental groups
    \[
    1 \to \pi_1(E_C) \to \pi_1(U_{0,g}) \to \Delta_g \to 1,
    \]
    where $\Delta_g := \pi_1^{orb}(\cH_g)$ is the hyperelliptic mapping class group. The asserted identification $\pi_1(E_C) \cong \Heis(H_1(C;\Z)) \times (\Z/2\Z)^2$ follows from \Cref{prop:Ubundle} in combination with \Cref{lemma:thetaheis}.
\end{proof}

\subsection{The principal bundle $\mathbf{E_C}$ in the case $\mathbf{n = g}$}\label{SS:ECng}

To prove \Cref{main:bundle} in the case $n = g$, we must understand $E_C/F_C$ as a principal bundle. Note that by \Cref{prop:Ubundle}, $F_C \cong C$. In \Cref{lemma:ECg}, we determine $\pi_1(E_C)$ as a central extension of $\pi_1(C)$, thereby finishing the proof of \Cref{main:bundle}. The proof of \Cref{lemma:ECg} is built on an analysis of the structure of the structure group $\Aut(\bF_g)$ (as $n = g \ne 0$, we drop the prime subscript of the previous section) - we must understand both its algebraic structure (\Cref{lemma:autFg}), as well as its topological structure (\Cref{lemma:Apres}).

\begin{lemma}\label{lemma:autFg}
    Let $\FAut(\bF_g)$ denote the subgroup of $\Aut(\bF_g)$ consisting of fiber-preserving automorphisms. Then there are short exact sequences
    \[
    1 \to \FAut(\bF_g) \to \Aut(\bF_g) \to \PGL_2(\C) \to 1
    \]
    and
    \[
    1 \to \C^{g+1} \to \Aut(\bF_g) \to (\GL_2(\C)/\mu_g) \to 1,
    \]
    where $\mu_g$ denotes the group of $g^{th}$ roots of unity (embedded in $\GL_2(\C)$ as a subgroup of scalar matrices).
    Moreover, $\FAut(\bF_g)$ admits a description as a semi-direct product
    \[
    \FAut(\bF_g) \cong H^0(\bP^1; \cO(g)) \semidirect H^0(\bP^1, \cO^*) \cong \C^{g+1} \semidirect \C^*.
    \]
\end{lemma}
\begin{proof}
    The second extension describing $\Aut(\bF_g)$ is given in \cite[Lemma 3.1]{JMG}. We will derive the other two assertions. As the fibers of the ruling of $\bF_g$ are the only irreducible curves of zero self-intersection, any automorphism of $\bF_g$ preserves the ruling and so descends to an automorphism of the base $\bP^1$. Taking this induced map induces the quotient $\Aut(\bF_g) \to \Aut(\bP^1) \cong \PGL_2(\C)$, and the kernel is the subgroup $\FAut(\bF_g)$ by definition. This gives the first extension above.

    To understand the structure of $\FAut(\bF_g)$, we start from the description $\bF_g = \bP(\cO \oplus \cO(-g))$. While this splitting is not canonical, the trivial sub-bundle $\cO$ is, and hence any fiberwise automorphism of the total space of the unprojectivized vector bundle must be of the form
    \[
    \begin{pmatrix}
    a & b\\ 0 & c
    \end{pmatrix}
    \]
    for $a, c \in H^0(\bP^1, \cO^*)$ and $b \in H^0(\bP^1, \Hom(\cO(-g),\cO)) \cong H^0(\bP^1, \cO(g))$. After projectivizing, this yields the isomorphism $\FAut(\bF_g) \cong H^0(\bP^1, \cO(g)) \semidirect H^0(\bP^1, \cO^*)$.
\end{proof}

\begin{lemma}\label{lemma:Apres}
    Setting $A \coloneq \pi_1(\Aut(\bF_g))$, there is a presentation of the form
    \[
    A \cong \Z\pair{v,w}/(v^g = w^2).
    \]
\end{lemma}
\begin{proof}
From the description in \Cref{lemma:autFg} of $\Aut(\bF_g)$ as an extension of $\GL_2(\C)/\mu_g$ by $\C^{g+1}$, we see that $\pi_1(\Aut(\bF_g)) \cong \pi_1(\GL_2(\C)/\mu_g)$. The covering space $\GL_2(\C) \to \GL_2(\C)/\mu_g$ is classified by the short exact sequence
\[
1 \to \pi_1(\GL_2(\C)) \to \pi_1(\GL_2(\C)/\mu_g) \to \mu_g \to 1.
\]
Both the kernel and the quotient are cyclic; let $w \in \pi_1(\GL_2(\C))$ be a generator, and let $v \in \pi_1(\GL_2(\C)/\mu_g)$ project to a generator of $\mu_g$. This shows that $A$ is a quotient of the free abelian group $\Z\pair{v,w}$ (recall that the fundamental group of any topological group is abelian). 

To derive the relation $v^g = w^2$, let $\tilde v:[0,1] \to \GL_2(\C)$ be the path given by $\tilde v(t) = e^{2 \pi i/g t}I$; such $\tilde v$ projects to a loop in $\GL_2(\C)/\mu_g$ that represents $v$. Its $g^{th}$ power is represented by the loop $e^{2 \pi i t} I$ (for $t \in [0,1]$). Now the isomorphism $\pi_1(\GL_2(\C)) \to \Z$ is realized at the space level by $\det: \GL_2(\C) \to \C^*$; this shows that $v^g$ represents the square of a generator, as claimed. It is easy to see that this is an exhaustive set of relations: any trivial $v^a w^b$ projects to a trivial element of $\mu_g$, so that $a = gk$. Then $v^g = w^2$ shows that $w^{b+2k} = 1$, but $w$ has infinite order.
\end{proof}

\begin{lemma}\label{lemma:ECg}
    In the case $n = g$, the space $E_C$ admits the following description of its fundamental group as a central extension:
    \[
    1 \to A \to \pi_1(E_C) \to \pi_1(C) \to 1,
    \]
    where $A = \pi_1(\Aut(\bF_g))$ is the abelian group $A \cong \Z\pair{v,w}/(v^g = w^2)$. The Euler class of this extension is $v^2 \in H^2(C;A) \cong A$. 
\end{lemma}
\begin{proof}
 As $E_C$ is a principal $\Aut(\bF_g)$-bundle over the aspherical space $C$ and (as seen in \Cref{lemma:autFg}) $\Aut(\bF_g)$ is connected, the fiber sequence for $E_C/C$ induces the short exact sequence
 \begin{equation}
     \label{eq:ECses}
     1 \to \pi_1(\Aut(\bF_g)) \to \pi_1(E_C) \to \pi_1(C) \to 1.
 \end{equation}
 According to \cite[Fact 1.4A]{nori}, under such circumstances, such an extension is necessarily central. 
 
 It remains to understand the extension class in $H^2(C;A) \cong A$. To do so, we consider the problem of specifying a local smooth section of $E_C/C$. Let $U \subset C$ be an open set. A section of $E_C$ over $U$ is given by a family of embeddings $i_p: C \to \bF_g$ (for $p \in U$), such that $i_p$ is induced by the line bundle $\cO_C(p)$ (in the setting of \Cref{lemma:hyptohirz}). Firstly, it is necessary to normalize such a family by specifying coordinates on the base $\bP^1$ under the projection $\pi: \bF_g \to \bP^1$. This can be accomplished by choosing a specific $h: C \to \bP^1$ such that $h = \pi \circ i_p$ for all $p$. Fix some such choice once and for all, and say that a (family of) embedding(s) $i_p$ is ``$h$-normalized'' if $\pi \circ i_p = h$.
 
  Recall from \Cref{lemma:hyptohirz} that $p$ can be recovered from an embedding $i: C \to \bF_g$ as the intersection $C \cdot S$ with the negative section. Thus, $i_p$ corresponds to an embedding $i_p^\circ: C\setminus\{p\} \to \bF_g \setminus S$. 
 For a line bundle $L$, let $E(L)$ denote the total space; then $\bF_g \setminus S \cong E(\cO_{\bP^1}(g))$ (see, e.g. \cite[Lemma 2.1]{fauntleroy}). Under this identification, $i_p^\circ$ acquires a simple pole at $p$. As $h^*(\cO_{\bP^1}(g)) \cong H^{\otimes g}$, it follows that the set of $h$-normalized embeddings above $p \in C$ is given by the set of genuinely meromorphic sections of $H^{\otimes g}(p)$.

To construct a smooth section of $E_C$ over some $U \subset C$, we must therefore construct a smoothly-varying family of genuinely meromorphic sections of $H^{\otimes g}(p)$ as $p \in U$ varies. To ensure a section of $H^{\otimes g}(p)$ is genuinely meromorphic, we examine the principal part. For a meromorphic section $\phi$ of a line bundle $L$ on $C$, the principal part at $p\in C$ corresponds to the functional on $L^* \otimes K_C$ that sends $\psi \otimes \omega$ to the residue of the form $\phi \otimes \psi \otimes \omega$ at $p$. Thus for varying $p \in C$, the bundle of principal parts is given by $H^{\otimes g}\otimes K_C^*$, and as $K_C \cong H^{\otimes (g-1)}$ for $C$ hyperelliptic, we conclude that the bundle of principal parts is isomorphic to $H$.

 Having understood the structure of $E_C/C$ as a principal bundle, we can now determine the Euler class of \eqref{eq:ECses} via sheaf cohomology. Let $\sE$ be the sheaf of smooth sections of the principal bundle $E_C/C$. Set $G = \Aut(\bF_g)$, and define $\sG$ as the sheaf of smooth $G$-valued functions on $C$. Then $\sE$ is a $\sG$-torsor, and so is classified by an element $[\sE]$ in the nonabelian sheaf cohomology $H^1(C,\sG)$ (a pointed set).

 The extension class of \eqref{eq:ECses} can be determined in this framework as follows. Let $\tilde G$ be the universal cover of $G$; then there is a short exact sequence of topological groups
 \[
 1 \to A \to \tilde G \to G \to 1. 
 \]
 Letting $\tilde \sG$ denote the sheaf of smooth $\tilde G$-valued functions, there is a short exact sequence of sheaves
 \[
 1 \to \underline{A} \to \tilde \sG \to \sG \to 1
 \]
 and hence a connecting map
 \[
 \delta: H^1(C, \sG) \to H^2(C, A);
 \]
unpacking the definitions shows that $\delta([\sE])$ is the extension class of \eqref{eq:ECses}.

To compute the extension class, we first note that since we are considering $h$-normalized embeddings, the problem of giving a smooth $\Aut(\bF_g)$-valued function on $C$ reduces to giving a smooth $\FAut(\bF_g)$-valued function. Via the description
\[
\FAut(\bF_g) \cong H^0(\bP^1, \cO(g)) \semidirect H^0(\bP^1, \cO^*) 
\]
 of \Cref{lemma:autFg}, this can in turn be formulated sheaf-theoretically. Define $\sR$ as the sheaf of {\em smooth} $\C$-valued functions on $C$, and likewise let $\sR^\times$ be the sheaf of smooth $\C^*$-valued functions. Define $\sF$ as the sheaf of smooth functions valued in $\FAut(\bF_g)$. Then there is a short exact sequence
 \[
 1 \to H^0(\P^1, \cO(g)) \otimes_\C \sR \to \sF \to \sR^\times \to 1.
 \]
 As $\sR$ is fine, it has vanishing $H^1$, and so any $\sF$-torsor is determined by its associated $\sR^\times$ torsor. But by construction, this torsor is nothing but the smooth line bundle on $C$ determined by the principal part of the meromorphic sections of $H^{\otimes g}(p)$. As discussed above, this is the line bundle $H$. We conclude that the extension class $\delta([\sE])$ can be identified with $c_1(H)$; it remains to determine the map $H^2(C; \Z) \to H^2(C; A)$ that induces such an identification.

 For this, we consider the diagram
\[
\begin{array}{ccccccccc}
0&\longrightarrow&\underline{\mathbb Z}&\longrightarrow&\mathscr R
&\xrightarrow{\exp(2\pi i\,\cdot)}&\mathscr R^\times&\longrightarrow&1\\
&&\downarrow j&&\downarrow\widetilde i&&\downarrow i\\
1&\longrightarrow&\underline A&\longrightarrow&\widetilde{\mathscr G}
&\longrightarrow&\mathscr G&\longrightarrow&1,
\end{array}
\]
where $i: \sR^\times \to \sG$ embeds $\C^*$ into the corresponding factor in $\FAut(F_g)$, which is in turn identified with the subgroup of scalar matrices in $\GL_2(\C)/\mu_g$. The map $\tilde i$ is the exponential lift. As seen in \Cref{lemma:Apres}, $v \in \pi_1(\Aut(\bF_g))$ is represented by a primitive loop of scalar matrices, and so $j: \Z \to A$ sends $1$ to $v$. As $c_1(H) = 2 \in H^2(C;\Z)$, it follows that the extension class $\delta([\sE]) \in H^2(C;A)$ is given as $2v$, as claimed.
\end{proof}

\begin{proof}[Proof of \Cref{main:bundle}, case $n = g$]
    The short exact sequence of orbifold fundamental groups
    \[
    1 \to \pi_1(E_C) \to \pi_1(U_{n,g}) \to \Delta_g \to 1
    \]
    can again be obtained from \Cref{prop:Ubundle}, while the precise description of $\pi_1(E_C)$ now follows from \Cref{lemma:ECg}.
\end{proof}

 \section{The smooth universal family}\label{section:universalfamily}
 
The analysis carried out in \Cref{lemma:hyptohirz} shows that for $1 \le d \le g+1$, every line bundle $L \in \Pic^d(C)$ determines an embedding $i: C \to \bF_n$, with the index $n$ determined by the behavior of $h^0(L\otimes H^r)$ as $r$ varies. On the other hand, the formula in \Cref{lemma:hyptohirz} shows that the index $n$ is well-defined mod $2$ over all of $\Pic^d(C)$, and hence the Hirzebruch surfaces are {\em diffeomorphic}. In this section, we show that these assemble into a smooth fiber bundle which is equipped with a fiberwise embedding of $C$, and examine some of the basic properties of this universal family.

\begin{lemma}\label{lemma:Mexists}
    Fix $C$ a hyperelliptic curve of genus $g \ge 2$, and a degree $1 \le d \le g+1$. Then there is a projective $g+2$-fold $\sM$ and a morphism $\pi: \sM \to \Pic^d(C)$ which as a map of smooth manifolds realizes $\sM$ as the total space of a smooth $\bF_n$-bundle for some $n \in \{0,1\}$ (specifically, $n \equiv g+1-d \pmod 2$). Moreover, there is an embedding $\sI: C\times \Pic^d(C) \to \sM$ restricting on each fiber of $\pi$ to an embedding of $C$.
\end{lemma}
\begin{proof}
    Let $h: C\to \bP^1$ be the hyperelliptic map. Let $\sL$ be a Poincar\'e line bundle on $C \times \Pic^d(C)$, i.e. a line bundle such that the pullback to $C \times \{L\}$ is isomorphic to $L$ for every $L \in \Pic^d(C)$ (for the construction of $\sL$, see \cite[p. 328]{GH} or \cite[IV.2]{ACGH}). For later use (\Cref{lemma:Mcohomology}), normalize $\sL$ as follows: let $p_0 \in C$ be a distinguished point, and normalize $\sL$ so that the restriction to $\{p_0\} \times \Pic^d(C)$ is trivial. Let 
    \begin{equation}
        \label{eq:Vdef}
        \sV:= (h \times \id)_*\sL
    \end{equation}
     be the rank-$2$ vector bundle on $\bP^1 \times \Pic^d(C)$ given by pushforward along the degree-$2$ map $h \times \id: C \times \Pic^d(C) \to \bP^1 \times \Pic^d(C)$. Define
    \[
    \sM\coloneq \bP(\sV).
    \]
    By construction, there is a projection $\sM \to \bP^1 \times \Pic^d(C) \to \Pic^d(C)$ which as a map of smooth manifolds is a proper surjective submersion, hence by Ehresmann's theorem yields a fiber bundle over $\Pic^d(C)$ with fibers diffeomorphic (though not all biholomorphic) to some fixed Hirzebruch surface $\bF_n$. Carrying out the construction of \Cref{lemma:hyptohirz} in this family, we obtain a family of embeddings $\mathscr{I}: C \times \Pic^d(C) \to \sM$ which is easily seen to assemble into a global embedding of $C \times \Pic^d(C)$.
\end{proof}

\begin{remark}\label{remark:normalized}
   The above construction gives a second projection $p: \sM \to \bP^1$, which restricts to a projection $p \circ \sI: C \times \Pic^d(C) \to \bP^1$. For later use (\Cref{lemma:F0epsbarkernel}), we note here that $p \circ \sI$ is {\em constant in the $\Pic^d(C)$-factor}, as indeed it is easy to see it coincides with the composition $C\times \Pic^d(C) \to C \to \bP^1$, with the first map being the projection and the second map the hyperelliptic map $h: C \to \bP^1$.
\end{remark}

As seen in \Cref{prop:Ubundle}, in translating from the ``intrinsic'' setting of a family of line bundles on $C$ to the ``explicit'' setting of a family of embeddings of some fixed hyperelliptic curve $C$ into $\bF_n$, one must be careful to specify the isomorphisms identifying the abstract ruled surfaces constructed in the intrinsic setting with a fixed model $\bF_n$. Having constructed the smooth universal family $\sM/\Pic^d(C)$, this is all the more important, as the isomorphism type of the Hirzebruch surface can vary. The result below shows how this is done in the case $d = g+1$ (where $C$ is embedded into $\bF_0$), in the case of a family that limits onto an embedding $C \into \bF_2$ at the central fiber; it will play a foundational role in the proof of \Cref{main:P1xP1}.

\begin{lemma}\label{lemma:holtriv}
    Let $\Delta \subset \C$ denote the closed unit disk, and let $\Delta^*\subset \Delta$ denote the complement of the origin. Let $j: \Delta \to \Pic^{g+1}(C)$ be a holomorphic map intersecting $\Theta'$ transversely at $t = 0$ in $\Delta$ and nowhere else. Let $M$ be the complex threefold produced by pulling the universal family $\sM \to \Pic^d(C)$ of \Cref{lemma:Mexists} back along $j$, and let $M^*$ denote the restriction of $M$ over $\Delta^*$. Then there is a distinguished family of biholomorphisms 
    \[
    \Phi: M^* \to \bF_0 \times \Delta^*
    \]
    which forms a torsor over the group $\Hol(\Delta, \PGL_2(\C))$ of holomorphic maps $\Delta \to \PGL_2(\C)$.
\end{lemma}
\begin{proof}
    Continuing in the setting of \Cref{lemma:Mexists}, we consider the vector bundle $(\id \times j)^*\sV$ on $\bP^1 \times \Delta$. For $t \in \Delta$, define $V_t$ as the restriction to $\bP^1 \times \{t\}$, so that by the argument of \Cref{lemma:hyptohirz}, $V_t \cong \cO \oplus \cO$ for $t \ne 0$ and $V_0 \cong \cO(1) \oplus \cO(-1)$.
    
    Consider the projection $\pi: \bP^1 \times \Delta \to \Delta$. We claim that when viewed as a locally free sheaf, $(\id \times j)^*\sV$ on $\bP^1 \times \Delta$ pushes forward under $\pi$ to a rank-$2$ vector bundle $W$ on $\Delta$. This follows since $h^0(\bP^1, V_t) = 2$ and $h^1(\bP^1, V_t) = 0$ for all $t \in \Delta$, including $t = 0$. 

    After possibly shrinking $\Delta$, there exists a holomorphic frame $e_1, e_2$ for $W$ defined over all of $\Delta$, corresponding to a basis of global sections $s_1, s_2 \in H^0(\bP^1 \times \Delta, (\id \times j)^* \sV)$. For $t \in \Delta^*$, such sections restrict to give a basis of global sections of $V_t \cong \cO \oplus \cO$; normalize $s_1$ and $s_2$ so that they are valued in the first (resp. second) factor. On a trivial bundle, such sections further restrict to a basis over every point in $\bP^1$, and hence give a distinguished identification 
    \[
    f_t: \bP(V_t) \cong \bF_0
    \]
    that varies holomorphically with $t$, as claimed. This construction depends only on the choice of holomorphic frame made above, and indeed only on the projective class of such a frame, showing that the space of choices forms a torsor over $\Hol(\Delta, \PGL_2(\C))$ as claimed.
\end{proof}

\begin{remark}\label{remark:mcgtriv}
    There is a crucial but counter-intuitive subtlety about these constructions. Any fiber bundle over the disk is trivial, and so with $M$ as in \Cref{lemma:holtriv}, there is a smooth trivialization $M \cong \bF_n \times \Delta$. However, in the case $d = g+1$ under consideration, this trivialization need not be compatible with a biholomorphism $\Phi: M^* \to \bF_0 \times \Delta^*$ as produced in \Cref{lemma:holtriv}. Under the smooth trivialization $M \cong \bF_0 \times \Delta$, restricting the embedding $i: C\times \Delta$ to $C \times \partial \Delta$ and projecting into $\bF_0$ yields a loop in $\Gr(C, \bF_0)$ which is null-homotopic, as it extends over $C \times \Delta$. {\em But this is not the loop given by the embedded monodromy $\epsilon_*$} - this is constructed by projecting onto $\bF_0$ via $\Phi$. 

    On the other hand, the existence of a smooth trivialization $M \cong \bF_0 \times \Delta$ {\em does} imply that the relative monodromy $\bar \epsilon_*$ of the loop around $\partial \Delta$, valued in $\pi_0(\Diff(\bF_0, C))$, is trivial. To see this, recall that a loop in $\Gr(C,\bF_0)$ determines a class in $\pi_0(\Diff(\bF_0,C))$ via isotopy extension, or equivalently, in the setting of \Cref{theorem:difffibr}, by applying the connecting map $\delta: \pi_1(\Gr(C,\bF_0)) \to \pi_0(\Diff(\bF_0,C))$ in the long exact sequence for the fibration $\Diff(\bF_0) \to \Gr(C,\bF_0)$. The loop in $\Gr(C,\bF_0)$ coming from $\Phi$ differs from the null-homotopic loop (coming from the smooth trivialization over $\Delta$) by the fiberwise application of a loop of diffeomorphisms of $\bF_0$ (relating the two trivializations of $M^*$). Such an element lies in the preceding term $\pi_1(\Diff(\bF_0))$ in the exact sequence, and so is ignored by $\delta$.
\end{remark}

We now consider the cohomology of $\sM$. We are particularly interested in the pullback $\sI^*: H^*(\sM;\Z) \to H^*(C \times \Pic^d(C);\Z)$, as this will be the primary mechanism by which we study the embedded monodromy of this family of curves. In preparation, we fix some notation that will be used throughout the paper. Assuming integral coefficients, define
\begin{equation}
\label{eq:etay}
    \eta \in H^2(C) \qquad \mbox{and} \qquad y \in H^2(\P^1)
\end{equation}
as the generators induced by the orientations compatible with the holomorphic structures on $C$ and $\P^1$. Next, recall the K\"unneth decomposition
\[
H^2(C \times \Pic^d(C)) \cong H^2(C) \oplus \left(H^1(C) \otimes H^1(\Pic^d(C))\right) \oplus H^2(\Pic^d(C)).
\]
Now there is a canonical identification $H^1(\Pic^d(C)) \cong H^1(C)$, and after applying the canonical isomorphism $H^1(C) \cong H_1(C)$ furnished by the intersection pairing, the middle term above can be canonically identified with $\End(H^1(C))$. In particular, there is a distinguished class 
\begin{equation}
    \label{eq:gamma}
    \gamma \in H^2(C \times \Pic^d(C))
\end{equation}

corresponding to the identity in $\End(H^1(C))$.

\begin{lemma}\label{lemma:Mcohomology}
    The cohomology algebra $H^*(\sM)$ is generated by $\pi^*(H^*(\bP^1 \times \Pic^d(C)))$ together with the class $\xi = c_1(\cO_{\sM}(1))$. Under the universal embedding $\sI: C \times \Pic^d(C) \to \sM$, 
    \[
    \sI^*(\xi) = d \eta + \gamma
    \]
 and 
    \[
    \sI^*(y) = 2 \eta.
    \]
    The pullback $\sI^*(H^*(\Pic^d(C)))$ coincides with the embedding $H^*(\Pic^d(C)) \to H^*(C \times \Pic^d(C))$ given by pullback along the projection.
\end{lemma}
\begin{proof}
    There is a commutative diagram
    \[
    \xymatrix{
    C \times \Pic^d(C) \ar[r]^{\sI} \ar[dr]_{h \times \id} & \sM \ar[d]^\pi\\
    & \P^1 \times \Pic^d(C),
    }
    \]
    showing in particular that $\sI^* \pi^*(H^*(\P^1 \times \Pic^d(C))) = (h \times \id)^*(H^*(\P^1 \times \Pic^d(C)))$, from which the latter two claims easily follow.

    The computation of $\sI^*(\xi) = \sI^*(c_1(\cO_\sM(1)))$ follows from the identification
    \[
    \sI^*(\cO_\sM(1)) \cong \sL,
    \]
    where $\sL$ is the normalized Poincar\'e line bundle over $C \times \Pic^d(C)$ previously encountered in \Cref{lemma:Mexists}. This identification is essentially tautological, as the embedding $\sI$ is constructed in the first place from the vector bundle $\sV = (h \times \id)_* \sL$.
     Consequently, $\sI^*(\xi) = c_1(\sL)$. 
     
     The Chern class $c_1(\sL)$ can be computed as follows. By construction, the restriction of $\sL$ to a fiber $(C,L)$ of $C \times \Pic^d(C)$ is the degree-$d$ line bundle $L$, from which it follows that the $H^2(C)$-component of $c_1(\sL)$ is $d \eta$. It is a basic property of Poincar\'e line bundles that the $H^1(C) \otimes H^1(\Pic^d(C))$-component of $c_1(\sL)$ is identified with $\id \in \End(H^1(C))$ - see \cite[p. 328]{GH}. Lastly, the $H^2(\Pic^d(C))$-component can be computed by restricting to the distinguished fiber $\{p_0\} \times \Pic^d(C)$, and this is trivial by our prior choice of normalization.
\end{proof}

We also have the following simple but quite useful fact. Here and throughout, for a compact oriented submanifold $Y \subset X$ of (real) codimension $k$, the notation 
    \[
    PD(Y) \in H^k(X)
    \]
    denotes the Poincar\'e dual.

\begin{lemma}\label{lemma:xiPDS}
    Let $L \in F_C \subset \Pic^d(C)$ be given, and let $i: \bP(h_*L) \to \sM$ denote the inclusion of the fiber. Then 
    \[
    i^*(\xi) = PD(S),
    \]
    where $\xi =  c_1(\cO_\sM(1))$ as above, and $PD(S)$ denotes the Poincar\'e dual to the distinguished section $S \subset \bP(h_*L)$.
\end{lemma}
\begin{proof}
    Identify $\bP(h_*L)\cong \bF_n$ for $n = g + 1 -d$. The pullback of $\cO_\sM(1)$ along $i$ is $\cO_{\bF_n}(1)$. The claim now follows from general principles. Let $p:V \to B$ be a rank-$2$ vector bundle, and let $s: B \to V$ be a nonvanishing section. Such $s$ determines a section $S$ of the projective bundle $\pi: \bP(V) \to B$ by fiberwise taking the quotient by the codimension-one subspace spanned by $s(t)$ (for $t \in B$). But via the universal quotient $\pi^*V \to \cO_{\P(V)}(1)$, such $s$ also determines a section of $\cO_{\P(V)}(1)$, whose vanishing locus is $S$ by construction. Applying this to the distinguished section $S$ of $\bF_n$ yields the result.
\end{proof}

Finally, we consider the Chern classes of $\sV$.

\begin{lemma}\label{lemma:chV}
    With $\sV$ defined as in \eqref{eq:Vdef}, the Chern classes of $\sV$ are computed in $H^*(\P^1 \times \Pic^d(C))$ as
    \[
    c_1(\sV) = (d-g-1)y, \qquad c_2(\sV) = y \theta,
    \]
    where $\theta \in H^2(\Pic^d(C))$ is the class of the principal polarization.
\end{lemma}
\begin{proof}
    This can be computed by the Grothendieck-Riemann-Roch formula, as applied to the line bundle $\sL$ on $C \times \Pic^d(C)$ and the morphism $h \times \id: C \times \Pic^d(C) \to \P^1 \times \Pic^d(C)$. As this morphism is finite, the higher direct images of $\sL$ vanish, and so the formula specializes to
    \[
    \ch(\sV)\td(\P^1 \times \Pic^d(C)) = (h \times \id)_*(\ch(\sL)\td(C \times \Pic^d(C))).
    \]
    As was seen in the proof of \Cref{lemma:Mcohomology}, $c_1(\sL) = d\eta + \gamma$. After unpacking the identification of $\gamma \in H^1(C) \otimes H^1(\Pic^d(C))$ with $\id \in \End(H_1(C))$, one computes 
    \[
    \gamma^2 = -2 \eta \theta \in H^2(C) \otimes H^2(\Pic^d(C)),
    \]
    where $\theta$ is the class of the principal polarization; we also have $\eta^2 = \eta \gamma = \gamma^3 = 0$. Thus $\ch(\sL) = e^{c_1(\sL)}$ is given as
    \[
    \ch(\sL) = 1 + d \eta + \gamma -\eta \theta.
    \]
    As $\Pic^d(C)$ is parallelizable, its contribution to the Todd classes is trivial, so that $\td(C \times \Pic^d(C))=\td(C) = 1 + (1-g)\eta$ and likewise $\td(\P^1\times \Pic^d(C)) = 1+ y$. 

    The pushforward $(h \times \id)_*: H^k(C \times \Pic^d(C)) \to H^k(\P^1 \times \Pic^d(C))$ is by definition given by conjugating the pushforward in homology by Poincar\'e duality. Thus
    \[
    (h \times \id)_*(1) = 2, \qquad (h \times \id)_*(\eta) = y, \qquad (h \times \id)_*(\gamma) = 0,
    \]
    with the latter holding since $\gamma \in H^1(C) \otimes H^1(\Pic^d(C))$ and hence 
    \[
    (h\times\id)_*(\gamma) \in H^1(\P^1) \otimes H^1(\Pic^d(C)) = 0. 
    \]
    Putting all of this together,
    \[
    \ch(\sV)(1+y) = 2 + (d+1-g)y -y \theta.
    \]
    Since $y^2 = 0$, multiplying by $(1-y)$ gives
    \[
    \ch(\sV) = 2 + c_1(\sV) + \tfrac{1}{2}c_1(\sV)^2 - c_2(\sV) = 2 + (d-g-1)y - y\theta,
    \]
    from which we see $c_1(\sV) = (d-g-1)y$ and hence $c_1(\sV)^2 = 0$, so that $c_2(\sV) = y\theta$.
\end{proof}

\section{\Cref{main:P1xP1}: hyperelliptic linear systems on $\bF_0$}\label{section:F0}

In this section we study the monodromy of the linear system $\abs{2S + (g+1)F}$ on $\bF_0 \cong \P^1 \times \P^1$, culminating in the proof of \Cref{main:P1xP1}. Recall the flavor of the results: we show that the relative monodromy
\[
\bar \epsilon_*: \pi_1(U_{0,g}) \to \pi_0(\Diff(\bF_0,C))
\]
is non-injective, but that the enhancement to the embedded monodromy
\[
\epsilon_*: \pi_1(U_{0,g}) \to \pi_1(\Gr(C,\bF_0))
\]
{\em is} sensitive enough to be faithful. We study the relative monodromy in \Cref{SS:relativep1p1} (\Cref{lemma:F0epsbarkernel}), and the embedded monodromy in \Cref{SS:embeddedtheta,SS:embeddedpgl}.

\subsection{Relative monodromy detects $\mathbf{H_1}$}\label{SS:relativep1p1}

\begin{lemma}\label{lemma:F0epsbarkernel}
    The relative monodromy
    \[
    \bar \epsilon_*: \pi_1(U_{0,g}) \to \pi_0(\Diff(\bF_0, C_{0,g}))
    \]
    has kernel $\Z \times (\Z/2\Z)^2$, realized geometrically as the subgroup generated by the theta meridian (the center of $\Heis(H_1(C;\Z))$) and the generators of $\pi_1(\Aut'(\bF_0))$.
\end{lemma}
\begin{proof}
      Recall from \Cref{prop:Ubundle} that the orbifold fiber over $C$ is isomorphic to $E_C \cong F_C \times \Aut'(\bF_0)$, where $F_C \cong \Pic^{g+1}(C) \setminus \Theta'$. We use this identification to lift the theta meridian in $\Pic^{g+1}(C) \setminus \Theta'$ into $E_C$ and hence into $U_{0,g}$.

    We first establish that $\ker(\bar \epsilon_*)$ contains the subgroup $\Z \times (\Z/2\Z)^2 \le \pi_1(U_{0,g})$ generated by the theta meridian and $\pi_1(\Aut'(\bF_0))$. That $\pi_1(\Aut'(\bF_0)) \le \ker \bar \epsilon_*$ is straightforward: indeed, for any pair of smooth manifolds $(X,C)$, any loop of diffeomorphisms of the ambient manifold $X$ induces a trivial class in $\pi_0(\Diff(X,C))$ (in the setting of \Cref{section:LES}, this is a consequence of the exactness of $\pi_1(\Diff(X)) \to \pi_1(\Gr(C,X)) \to \pi_0(\Diff(X,C))$). 
    Letting $z$ denote the class of the theta meridian, $z$ is represented by a loop $\zeta: S^1 \to \Pic^{g+1}(C) \setminus \Theta'$, where $\Theta'$ is identified with $W_{g-1} \subset \Pic^{g-1}(C)$. Letting $\Delta \subset \C$ denote the unit disk as above, such $\zeta$ extends to a map $Z: \Delta \to \Pic^{g+1}(C)$, and so by \Cref{remark:mcgtriv}, $\bar \epsilon_*(z) = 1$ as claimed. 

    To establish the other containment, first observe that $\ker \bar \epsilon_* \le \ker \rho_*$, where $\rho_*: \pi_1(U_{0,g}) \to \Mod(\Sigma_g)$ is the intrinsic monodromy. By \Cref{main:bundle}, $\ker \rho_* \cong \Heis(H_1(C;\Z)) \times (\Z/2\Z)^2$, and by the previous paragraph, the restriction of $\bar \epsilon_*$ to this subgroup factors through $\Heis(H_1(C;\Z))/ \pair{z} \cong H_1(C;\Z)$. It thus suffices to show that the induced map $H_1(C;\Z) \to \pi_0(\Diff(\bF_0,C))$ is injective.

    To establish this, we consider the action of $\pi_0(\Diff(\bF_0,C))$ on the relative homology $H_2(\bF_0, C; \Z)$ (we suppress integer coefficients going forward). The long exact sequence for the pair $(\bF_0, C)$ contains the segment
    \[
    0 \to H_2(C) \to H_2(\bF_0) \to H_2(\bF_0, C) \to H_1(C) \to 0,
    \]
    realizing $H_2(\bF_0,C)$ as an extension of $H_1(C)$ by $H_2(\bF_0)/H_2(C)$. As $H_2(C)$ includes in $H_2(\bF_0)$ as the cyclic subgroup spanned by $2S + (g+1)F$, we see that the kernel of this extension is isomorphic to $\Z$ if $g$ is even, and to $\Z \oplus \Z/2\Z$ if $g$ is odd.
    
    The action on the quotient $H_1(C)$ is the usual symplectic representation of $\Mod(C)$, which is trivial in our setting. The action on the preceding term $H_2(\bF_0)$ is likewise trivial, as our diffeomorphisms arise from isotopy extension and hence are isotopic to the identity if $C$ is free to move. The action of $\bar \epsilon_*(H_1(C)) \le \pi_0(\Diff(\bF_0,C))$ on $H_2(\bF_0,C)$ is therefore valued in the subgroup of transvections $\Hom(H_1(C), H_2(\bF_0)/H_2(C)) \le \Aut(H_2(\bF_0,C))$. Let
    \[
    \tau: H_1(C) \to \Hom(H_1(C), H_2(\bF_0)/H_2(C))
    \]
    be the homomorphism induced by restricting $\bar \epsilon_*$ to $\Heis(H_1(C)) \le \pi_1(U_{0,g})$ and then descending to the quotient $H_1(C) \cong \Heis(H_1(C))/\pair{z}$. To conclude the argument, it suffices to establish the formula
    \begin{equation}
        \label{eq:tau}
         \tau(v)(w) = \pair{v,w}[F]
    \end{equation}
    (where $\tau(v)(w)$ is valued in $H_2(\bF_0)/H_2(C)$).

    To establish this, let $f \in \Diff(\bF_0,C)$ act trivially on $H_1(C)$ and $H_2(\bF_0)$. The image of $f$ in $\Hom(H_1(C), H_2(\bF_0)/H_2(C))$ can be given by 
    \[
    w \mapsto f_*\tilde w - \tilde w,
    \]
    where $\tilde w \in H_2(\bF_0,C)$ lifts $w \in H_1(C)$. 
    In our setting it is possible to construct an explicit representative of this cycle (this is reminiscent of Schnell's {\em tube mapping} \cite{schnell}, or of the construction of the Johnson homomorphism via the intersection form on the mapping torus \cite{johnsonsurvey}). Let $v \in H_1(C)$ be given, represented by a map $v: S^1 \to \Pic^{g+1}(C)$, giving rise (by \Cref{lemma:Mexists}) to an $S^1$-family of embeddings $i_{v}: C \times S^1 \to \bF_0 \times S^1$. Let  $f_v \in \Diff(\bF_0,C)$ be a lift of $\bar \epsilon_*(v)$ that restricts to the identity on $C$. 
    
    Now choose $\omega \subset C$ to be a curve representing $w \in H_1(C)$, and let $\tilde \omega \subset \bF_0$ be a surface in $\bF_0$ bounding $\omega \subset C$. The absolute cycle $f_v(\tilde \omega) - \tilde \omega$ representing $f_{v,*} \tilde w - \tilde w \in H_2(\bF_0)$ is then seen to be the image of the torus $(p_1 \circ i_v)(T)$
    obtained by pushing $\omega \subset C$ along the loop of embeddings $i_v$ (here, $p_1: \bF_0 \times S^1 \to \bF_0$ is the projection onto the first factor). Let
    \[
    T_{v,\omega} = (p_1 \circ i_v)_*[T]
    \]
    denote this cycle.
    
    Following \Cref{remark:normalized}, the projection of $C \times \{t\}$ onto $\bP^1$ along the $F$-ruling is constant in $t \in S^1$. Thus the restriction to $T_{v,\omega}$ has zero degree, showing that in $H_2(\bF_0) = \Z\pair{S,F}$,
    \[
    T_{v, \omega} = n_{v,\omega}[F]
    \]
    for some integer $n_{v,\omega}$.
    To establish \eqref{eq:tau}, we must show that $n_{v,\omega} = \pair{v,w}$. To see this, note that $n_{v,\omega}$ can be computed intersection-theoretically as 
    \[
    n_{v,\omega} = T_{v,\omega}\cdot S.
    \] 
    Recall from \Cref{lemma:xiPDS} that $S$ is Poincar\'e dual to the Chern class of $\cO_{\bF_0}(1)$, i.e. the pullback of $\xi \in H^2(\sM)$ to $\bF_0$, so that
    \[
    n_{v,\omega} = \pair{(p_1 \circ i_v)^*(c_1(\cO_{\bF_0}(1))), [T]}
    \] 
    
    By \Cref{lemma:Mcohomology}, $(p_1 \circ i_v)^*(c_1(\cO_{\bF_0}(1)))$ is given by the pullback of $(g+1)\eta + \gamma \in H^2(C \times \Pic^{g+1}(C))$ along the composition $T \to C \times S^1 \xrightarrow{\id \times v} C \times \Pic^{g+1}(C)$, where $\gamma \in H^1(C) \otimes H^1(\Pic^{g+1}(C))$ is identified with $\id \in \End(H_1(C))$ (recall \eqref{eq:gamma}). Thus 
    \[
    (\id \times v)^*(\gamma) = PD(v) \otimes \sigma \in H^1(C) \otimes H^1(S^1),
    \]
    where $\sigma \in H^1(S^1)$ is a generator.
    Pulling back along the inclusion $T = \omega \times S^1 \into C \times S^1$ sends $PD(v) \otimes \sigma$ to $\pair{v,w} \sigma_1 \otimes \sigma_2$, where $\sigma_i$ represents a generator of the $i^{th}$ factor $S^1 \subset T \cong S^1\times S^1$. This shows $n_{v,\omega} = \pair{v,w}$ as claimed.
    \end{proof}

\subsection{Embedded monodromy detects the theta meridian}\label{SS:embeddedtheta}
In this section and the next we consider the more sensitive {\em embedded monodromy}, the map $\epsilon_*: \pi_1(U_{0,g}) \to \pi_1(\Gr(C,\bF_0))$. We show that $\epsilon_*$ detects the kernel of $\bar \epsilon_*$ computed in \Cref{lemma:F0epsbarkernel}. In this section, we specialize to a study of the theta meridian; the kernel terms arising from $\Aut'(\bF_0)$ are studied in \Cref{SS:embeddedpgl}.

\begin{lemma}\label{lemma:thetameridianinforder}
    Let $z \in \pi_1(\Pic^{g+1}(C) \setminus \Theta')$ denote the theta meridian, lifted via \Cref{prop:Ubundle} to an element of $\pi_1(E_C) \le \pi_1(U_{0,g})$. Then under the embedded monodromy
    \[
    \epsilon_*: \pi_1(U_{0,g}) \to \pi_1(\Gr(C,\bF_0)),
    \]
    $\epsilon_*(z) \in \pi_1(\Gr(C,\bF_0))$ is an element of infinite order.
\end{lemma}

\begin{proof}
    By \Cref{lemma:F0epsbarkernel}, $\bar \epsilon_*(z) = 1$ and so {\em a fortiori} also $\rho_*(z) = 1$. Thus by \Cref{lemma:kerrholifts}, $\epsilon_*(z) \in \pi_1(\Gr(C,\bF_0))$ lifts to $\tilde \epsilon_*(z) \in \pi_1(\Emb(C, \bF_0))$. By (the proof of) \Cref{lemma:kerrholifts}, the natural map $\pi_1(\Emb(C,\bF_0)) \to \pi_1(\Gr(C,\bF_0))$ is injective, so that if $\tilde \epsilon_*(z) \in \pi_1(\Emb(C,\bF_0))$ has infinite order, then so does $\epsilon_*(z) \in \pi_1(\Gr(C,\bF_0))$. 
    
    We will show that $\tilde \epsilon_*(z)$ has infinite order in four steps. Most of the work will take place in the context of showing that $\tilde \epsilon_*(z) \ne 1$; that this moreover has infinite order will be an easy corollary.
    
    \para{Step 1: Setup} 
    The theta meridian $z: S^1 \to \Pic^{g+1}(C) \setminus \Theta'$ can be chosen so as to extend to a holomorphic map $Z: \Delta \to \Pic^{g+1}(C)$ that meets $\Theta'$ transversely at $t =0$. By \Cref{lemma:Mexists}, this produces a complex manifold $M$ that as a smooth manifold is an $S^2 \times S^2$-bundle $p: M \to \Delta$, together with a holomorphic map $i: C \times \Delta \to M$ that is a fiberwise embedding.

    We show that $\tilde \epsilon_*(z) \in \pi_1(\Emb(C, \bF_0))$ has infinite order in a slightly indirect way. By \Cref{lemma:holtriv}, choose a biholomorphism $\Phi: M^* \to \bF_0 \times \Delta^*$. Letting $p_0: \bF_0 \times \Delta^* \to \bF_0$ denote the projection, the composition 
    \[
    p_0 \circ \Phi \circ i: C \times S^1 \to \bF_0
    \]
    then yields a map $j: S^1 \to \Emb(C, \bF_0)$ representing $\tilde \epsilon_*(\tilde z)$ for some lift $\tilde z$ of $z$ from $\pi_1(F_C)$ to $\pi_1(E_C)$. Such $\tilde z$ need not be the lift of $z$ given by the splitting of \Cref{prop:Ubundle}- these may differ by an element of $(\Z/2\Z)^2$. But then their squares are necessarily equal, and so one has infinite order if and only if the other does. We therefore proceed by showing that $\tilde \epsilon_*(\tilde z)$ has infinite order.

    \para{Step 2: $\mathbf{N}$ and its cohomology}  Let $N$ be the holomorphic $3$-fold obtained by extending the trivial family $M^* \cong \bF_0 \times \Delta^*$ from $\Delta^*$ to the whole of $\bP^1\setminus \{0\}$. Then $N$ is a holomorphic family of ruled surfaces over $\bP^1_t$ with fiber $\bF_2$ over $t=0$ and fiber $\bF_0$ for $t \ne 0$ (momentarily there will be two distinct $\P^1$'s in play; we distinguish them with a subscript indicating their corresponding coordinate). By construction, 
    \[
    N = \bP(\hat \sV),
    \]
    where $\hat \sV$ is a rank-$2$ vector bundle over $\bP^1_y \times \bP^1_t$ extending (in the notation of \Cref{lemma:Mexists}) the rank-$2$ vector bundle $(\id \times Z)^*(\sV)$ over $\bP^1_y \times \Delta$. For $t \in \bP^1_t$, let $V_t$ be the vector bundle on $\bP^1_y \times \{t\}$ obtained by restricting $\hat \sV$. Let $\pi: N \to \bP^1_y \times \bP^1_t$ be the projection. In \Cref{lemma:chernsV,lemma:Ncohomology,lemma:p1TN}, we establish some lemmas concerning characteristic classes of various vector bundles on $N$.

    \begin{lemma}\label{lemma:chernsV}
        The Chern classes of $\hat \sV$ satisfy
        \[
        c_1(\hat \sV) = r\, t, \qquad c_2(\hat \sV) = s\, yt,
        \]
        where $r,s$ are integers with $s \ne 0$, and $y,t \in H^2(\bP^1_y \times \bP^1_t;\Z)$ denote the pullbacks of positive generators of $H^2(\bP^1_y;\Z)$ and $H^2(\bP^1_t;\Z)$, respectively. 
    \end{lemma}
    \begin{proof}
        Write $c_1(\hat \sV) = q\,y + r\,t$ for integers $q, r$. The integer $q$ can be shown to vanish by restricting $\hat \sV$ to any slice $\bP^1_y \times \{t\}$, yielding a vector bundle $V_t$ of degree zero (as $V_t \cong \cO_{\bP^1} \oplus \cO_{\bP^1}$ for $t \ne 0$ and $V_0 \cong \cO_{\bP^1}(1) \oplus \cO_{\bP^1}(-1)$).
        
        To see that $c_2(\hat \sV) \ne 0$, we compute the vanishing locus of a section explicitly. Returning to the setting of \Cref{lemma:holtriv}, let $\hat s_1 \in H^0(\bP^1_y \times \bP^1_t, \hat \sV)$ be the global section that restricts to the section $(1,0) \in V_t \cong \cO_{\bP^1} \oplus \cO_{\bP^1}$ for $t \in \bP^1_t \setminus \Delta$. Such $\hat s_1$ is nonvanishing for $t \ne 0$, since it is still part of a frame for the trivial bundle fiber $V_t$. 
        
        At $t = 0$, it restricts to a section of $V_0 \cong \cO_{\bP^1}(1) \oplus \cO_{\bP^1}(-1)$ and so is necessarily of the form $(s, 0)$ for $s$ a nonzero section of $\cO_{\bP^1}(1)$. Such $s$ has exactly one zero at some $y_0 \in \bP^1_y$, so the vanishing locus of $\hat s_1$ is the single point $(y_0, 0) \in \bP^1_y \times \bP^1_t$. An isolated zero necessarily has positive multiplicity, showing that the Chern number is nonvanishing as claimed.
    \end{proof}

 \begin{remark}
        An alternative approach to \Cref{lemma:chernsV} is to extend the universal computation in \Cref{lemma:chV} to the relative setting, and compute the Chern classes relative to the canonical trivializations of $\sM$ and related bundles that exist away from the theta locus. This makes it easy to obtain the Chern classes of $\hat \sV$, as we obtain $N$ by capping off a subfamily of $\sM$ over a disk along this preferred trivialization. While this may be conceptually and technically cleaner, we have chosen to take the more explicit approach presented here so as to avoid introducing even more algebro-geometric machinery than what is already present in what is supposed to be a topology paper.
    \end{remark}

    \begin{lemma}\label{lemma:Ncohomology}
    Let $y,t \in H^2(N;\Q)$ be the pullbacks of positive generators of $H^2(\bP^1_y;\Q)$ and $H^2(\bP^1_t;\Q)$, respectively, let $\xi_N = c_1(\cO_N(1))$, and define $x = \xi_N - \tfrac{r}{2}t$, where $r, s \in \Z$ are as in \Cref{lemma:chernsV}. Then 
    there is a ring isomorphism
    \[
    H^*(N;\Q) \cong \Q[x,y,t]/(y^2 = t^2 = 0, x^2 = -s\, yt).
    \]    
    \end{lemma}
    \begin{proof}
        In general, suppose $\mathscr{W}$ is a vector bundle of rank $2$ over a base $B$, and let $E = \bP(\mathscr{W})$; then $H^*(E;\Q)$ is generated by $c_1(\cO_E(1))$ as an algebra over $H^*(B;\Q)$. This shows that $H^*(N;\Q)$ is generated as a ring by $x,y,t$ and that $y^2 = t^2 =0$; it remains to obtain the relation $x^2 = -s\, yt$. 

        For this, we consider the tautological decomposition of $\pi^*(\hat \sV)$
        \begin{equation}\label{eq:tautological}
            0 \to K \to \pi^*(\hat \sV) \to \cO_N(1) \to 0,
        \end{equation}
        defining $K$ as the kernel of the tautological quotient $\pi^*(\hat \sV) \to \cO_N(1)$. Applying the Whitney sum formula and \Cref{lemma:chernsV} shows that
        \begin{equation}\label{eq:c1K}
            c_1(K) = r\, t - \xi_N 
        \end{equation}
    
        and
        \[
        c_2(\pi^*(\hat \sV)) = c_1(K) c_1(\cO_N(1)) = (r\,t-\xi_N)\xi_N.
        \]
        On the other hand, by \Cref{lemma:chernsV}, 
        \[
        c_2(\pi^*(\hat \sV)) = \pi^*(c_2(\hat \sV)) = \pi^*(s\,yt).
        \]
        Equating these descriptions and rearranging gives
        \[
        \xi_N^2 -r\,\xi_N t + s\,yt = 0.
        \]
        Defining $x\coloneq \xi_N-\tfrac{r}{2}t$ and factoring, the above expression transforms into $x^2 = -s\,yt$.
    \end{proof}

        \begin{lemma}\label{lemma:p1TN}
        The Pontryagin class of $TN$ is given by the formula
        \[
        p_1(TN) = -4s\, yt.
        \]
    \end{lemma}
    \begin{proof}
        From the construction $N = \bP(\hat \sV)$ as the total space of a $\bP^1$-bundle over $\bP^1_y \times \bP^1_t$, 
        \[
        TN \cong \pi^*(T(\bP^1_y \times \bP^1_t)) \oplus \omega,
        \]
        where $\omega$ is the vertical tangent bundle. From the tautological decomposition \eqref{eq:tautological},
        \[
        \omega \cong \Hom(K, \cO_N(1)),
        \]
        showing by \eqref{eq:c1K} and the definition $x \coloneq \xi_N - \tfrac{r}{2}t$ that
        \[
        c_1(\omega) = \xi_N - c_1(K) = \xi_N - (r\,t - \xi_N) = 2\xi_N -r\,t = 2x.
        \]
        Thus
        \begin{align*}
            p_1(TN) &= p_1(\pi^*(T(\bP^1_y \times \bP^1_t))) + p_1(\omega)\\
            &= 0 + c_1^2(\omega)\\
            &= 4x^2\\
            &= -4s\,yt,
        \end{align*}
      with the last equality holding by \Cref{lemma:Ncohomology}.
    \end{proof}

    \para{Step 3: A hypothetical null-homotopy} We now suppose for the sake of contradiction that the loop $\tilde \epsilon_*(\tilde z)$ is null-homotopic, so that the above map $j$ admits an extension $J: \Delta \to \Emb(C, \bF_0)$. Then $i: C \times \Delta \to M$ admits an extension $I: C \times \bP^1_t \to N$, determining a smoothly-embedded submanifold $I(C \times \bP^1_t) \subset N$ of (real) codimension two. We will show that such a submanifold cannot exist, by examining how the vector bundles studied in Step 2 above would behave when restricted to $I(C \times \bP^1_t)$. 

    To ease the notational load, we will use the same symbol ($\eta, y,t$, etc.) to denote a cohomology class and its pullback along an inclusion or projection map; we hope this does not lead to undue confusion.

    \begin{lemma}\label{lemma:ipullbackformula}
        Let $\eta, t \in H^2(C \times \bP^1_t;\Q)$ denote the pullbacks of positive generators of $H^2(C;\Q)$ and $H^2(\bP^1_t;\Q)$, respectively. Then under $I: C \times \bP^1_t \to N$, 
        \[
        I^*(x) = m \eta - \tfrac{s}{m}t, \quad I^*(y) = 2 \eta, \quad I^*(t) =t,
        \]
        for some $m \in \Q$; consequently
        \[
        I^*(xy) =-\tfrac{2s}{m} \eta t, \quad I^*(xt) = m \eta t, \quad I^*(yt) = 2 \eta t.
        \]
    \end{lemma}
    \begin{proof}
        The formula $I^*(t) = t$ is immediate. To compute $I^*(y)$, write $I^*(y) = b\eta + ct$. Restricting further to a fiber $C \times \{t\}$ shows that $b = 2$. Then pulling back the relation $y^2 = 0$ shows that $c = 0$, as desired.
        
        To compute $I^*(x)$, write $I^*(x) = m \eta + n t$ for some $m,n \in\Q$. The relation $x^2 = -s\, yt$ in $H^*(N;\Q)$ pulls back to show
        \[
        2mn\, \eta t = -2s\, \eta t.
        \]
       As $s \ne 0$ by \Cref{lemma:chernsV}, it follows that $m \ne 0$ and $n = -\tfrac{s}{m}$ as claimed. The formulas for the classes of degree $4$ now follow from \Cref{lemma:Ncohomology}. 
    \end{proof}

    \begin{lemma} \label{lemma:deltaformula}
        Define $\delta \in H^2(N;\Q)$ as the Poincar\'e dual to $I(C \times \bP^1_t)$. Then
        \[
        \delta = 2x + my - \tfrac{2s}{m}t,
        \]
        where $m$ is defined in \Cref{lemma:ipullbackformula}.
    \end{lemma}
    \begin{proof}
        Write
        \[
        \delta = a\,x + b\,y + c\,t
        \]
        for $a,b,c \in \Q$. Using the characteristic property of Poincar\'e duals and the formulas of \Cref{lemma:ipullbackformula},
        \begin{align*}
            &a = \pair{\delta y t,[N]} = \pair{I^*(yt), [C \times \bP^1_t] } = \pair{2 \eta t,[C \times \bP^1_t]} = 2,\\
            &b = \pair{\delta x t ,[N]} = \pair{I^*(xt),[C \times \bP^1_t]} = \pair{m \eta t, [C \times \bP^1_t]} = m,\\
            &c = \pair{\delta x y, [N]} = \pair{I^*(xy), [C \times \bP^1_t]} =  \pair{- \tfrac{2s}{m}\eta t,[{C \times \bP^1_t}]} = - \tfrac{2s}{m}.
        \end{align*}
    \end{proof}

\para{Step 4: Deriving the contradiction} We can now derive a contradiction to the existence of the embedding $I: C \times \bP^1_t \to N$, by computing $p_1(I^*(TN))$ and $I^*(p_1(TN))$ and seeing they cannot agree. On the one hand, the bundle $I^*TN$ on $C \times \bP^1_t$ splits as
\[
I^*TN \cong T(C \times \bP^1_t) \oplus \nu,
\]
where $\nu$ is the normal bundle. By the Whitney formula, the vanishing of the Pontryagin class on $C \times \bP_t^1$ and the equality $p_1(\nu) = e(\nu)^2$,
\[
p_1(I^*TN) = p_1(T(C \times \bP^1_t)) + p_1(\nu) = p_1(\nu) = e(\nu)^2.
\]
By Poincar\'e duality, $e(\nu) = I^*(\delta)$, so that by \Cref{lemma:Ncohomology,lemma:ipullbackformula,lemma:deltaformula,},
\begin{align*}
    p_1(I^*TN) &= e(\nu)^2 = I^*(\delta)^2 = I^*(2x + my - \tfrac{2s}{m}t)^2\\
        &= I^*(4m\, xy - \tfrac{8s}{m}\,xt -8s\,yt)\\
        &= -32s\,\eta t.
\end{align*}

On the other hand, by \Cref{lemma:p1TN} and \Cref{lemma:ipullbackformula},
\begin{align*}
I^*(p_1(TN)) &= I^*(-4s\, yt)\\
&= -8s\, \eta t.
\end{align*}
Recall from \Cref{lemma:chernsV} that $s \ne 0$; this shows the inequality $p_1(I^*(TN)) \ne I^*(p_1(TN))$, a contradiction.

This shows that $\tilde \epsilon_*(\tilde z) \in \pi_1(\Emb(C, \bF_0))$ is nontrivial. To show it has infinite order, observe that for any $k \ge 1$, the loop $\tilde \epsilon_*(\tilde z^k)$ can be represented by pulling back the original family over $\Delta$ along the $k^{th}$ power map $p_k: \Delta \to \Delta$. Under such a map, the Chern classes of the pullback of $\hat \sV$ scale by a factor of $k$, and in particular the second Chern class remains nonvanishing; the argument then proceeds {\em mutatis mutandis}.
\end{proof}

\subsection{Embedded monodromy detects the $\Z/2$ factors}\label{SS:embeddedpgl}
Finally here we show that $\epsilon_*$ also detects the subgroup $(\Z/2\Z)^2 \le \ker \bar \epsilon_*$ arising from $\pi_1(\PGL_2(\C))$. This is a consequence of a more general result.

\begin{lemma}\label{lemma:mapspace}
    Let $C$ be a closed oriented surface, let $g: C \to \bP^1$ be a continuous map of nonzero degree, and let $A: S^1 \to \PGL_2(\C)$ induce the nontrivial element of $\pi_1(\PGL_2(\C))$. Then the composition
    \[
    t \mapsto [A_t \circ g : C \to \bP^1]
    \]
    determines a nontrivial element in $\pi_1(\Map(C,\bP^1))$.
\end{lemma}
\begin{proof}
    Define 
    \[
    G: S^1 \times C \to S^1 \times \bP^1
    \]
    by $G(t,\cdot)= (t,A_t \circ g)$. Suppose for the sake of contradiction that $G$ admitted a continuous extension over the disk $\tilde G: \Delta \times C \to \Delta \times \bP^1$. Let $M$ be the $4$-manifold constructed by gluing two copies of $\Delta \times \bP^1$ along their common boundary, identifying $(t,z)$ in one copy with $(t,A_t^{-1}(z))$ in the other. 
    
    Topologically, $M$ is the nontrivial $S^2$-bundle over $S^2$, so there is a diffeomorphism $M \cong \bF_1$; in particular their cohomology rings are isomorphic. We let $PD(S), PD(F) \in H^2(M;\Z)$ denote the Poincar\'e duals of the cycles $S, F \subset \bF_1$, respectively. These satisfy the relations $PD(S)^2 = -PD(S) PD(F)$ and $PD(F)^2 = 0$. 

    Under the gluing used to construct $M$, the map $G: S^1\times C \to S^1 \times \bP^1$ appears constant on the boundary of the second copy of $\Delta \times \bP^1$, and so there is a further extension $\hat G: S^2 \times C \to M$. 
    We will obtain a contradiction by examining the map $\hat G^*: H^*(M;\Z) \to H^*(S^2 \times C;\Z)$. As usual, let $t, \eta\in H^2(S^2 \times \C;\Z)$ denote the pullbacks of positive generators of $H^2(S^2;\Z)$ and $H^2(C;\Z)$, respectively. 

    Defining $d \coloneq \deg(g)$, 
    \[
    \hat G^*(PD(S)) = d \eta + a t
    \]
    for some $a \in \Z$, and
    \[
    \hat G^*(PD(F)) = t,
    \]
    the latter holding by commutativity of the diagram
    \[
    \xymatrix{
    S^2 \times C \ar[r]^{\hat G} \ar[d] & M \ar[d]\\
    S^2 \ar@{=}[r] &S^2
    }
    \]
    and the definition of the classes $t \in H^2(S^2 \times C;\Z)$ and $PD(F) \in H^2(M;\Z)$.
    Pulling back the relation $PD(S)^2 = -PD(S) PD(F)$ then shows
    \[
    2ad = -d
    \]
    but as $d \ne 0$ by hypothesis and $a$ is integral, this is impossible.
\end{proof}

\begin{lemma}\label{lemma:pglinforder}
    The restriction of $\epsilon_*: \pi_1(U_{0,g}) \to \pi_1(\Gr(C_{0,g}, \bF_0))$ to the subgroup $\pi_1(\Aut'(\bF_0)) \cong (\Z/2\Z)^2$ is an injection.
\end{lemma}
\begin{proof}
As in \Cref{lemma:thetameridianinforder}, we lift $\epsilon_*: \pi_1(\Aut'(\bF_0)) \to \pi_1(\Gr(C,\bF_0))$ up to $\tilde \epsilon_*: \pi_1(\Aut'(\bF_0)) \to \pi_1(\Emb(C,\bF_0))$ (and as discussed there, it suffices to show that this lift is injective). Let $A: S^1 \to \Aut'(\bF_0)$ determine a nontrivial element. Then $A$ remains nontrivial under at least one of the two projections $\Aut'(\bF_0) \to \PGL_2(\C)$. The degrees of the corresponding maps $C \subset \bF_0 \to \bP^1$ are $2$ and $g+1$, both nonzero. The projection induces a continuous map $\Emb(C, \bF_0) \to \Map(C, \bP^1)$, and the image of $A$ under this map is then seen to determine a nontrivial element of $\pi_1(\Map(C,\bP^1))$ by \Cref{lemma:mapspace}.    
\end{proof}

\begin{proof}[Proof of \Cref{main:P1xP1}]
    The statement of \Cref{main:P1xP1} is an amalgamation of the results contained in \Cref{lemma:F0epsbarkernel,lemma:thetameridianinforder,lemma:pglinforder}.
\end{proof}

\section{\Cref{main:Fg}: hyperelliptic linear systems on $\bF_g$}\label{section:Fg}

We now turn to \Cref{main:Fg}, and consider the linear system $\abs{2S + (2g+1)F}$ on $\bF_g$. 

\subsection{A smooth extension of the embedded monodromy}
As seen in \Cref{lemma:Mexists}, the family of embeddings $C \to \bF_g$ parametrized by $C \subset \Pic^1(C)$ extends over the whole of $\Pic^1(C)$, at the cost of embedding $C$ in different Hirzebruch surfaces $\bF_n$ (though $n \equiv g \pmod 2$). At first glance, this might suggest that there is a map $C \to \Emb(C,\bF_g)$ that admits an extension to $\Pic^1(C) \to \Emb(C,\bF_g)$, and hence that the embedded monodromy $\tilde \epsilon_*: \pi_1(C) \to \pi_1(\Emb(C,\bF_g))$ factors through the abelian group $\pi_1(\Pic^1(C)) \cong H_1(C;\Z)$. 

This intuition is not precisely correct - indeed \Cref{lemma:inontriv} shows that $\tilde \epsilon_*$ does detect certain commutators. It can be repaired by considering the same mechanism underlying the nontriviality of the theta meridian in the previous section: an arbitrary $L \in \Pic^1(C)$ only furnishes an embedding into an abstract ruled surface $\bP(h_*L)$ and not a fixed model $\bF_n$ (note of course that the index $n$ depends on $L$, though its parity does not). For one thing, as seen in \Cref{lemma:ECg}, the principal bundle $E_C/C$ is nontrivial, so that it is not even possible to globally identify the family of abstract ruled surfaces over $C$ with a fixed model $\bF_g$ by biholomorphisms. Thus there is no map $C \to \Emb(C, \bF_g)$ hypothesized above, only a map $E_C \to \Emb(C, \bF_g)$. In extending to embeddings defined using general line bundles in $\Pic^1(C)$, there are no holomorphic identifications of the resulting ruled surface with $\bF_g$ at all; instead we must pass to the smooth category. This leads us to the following definition.

\begin{definition}
    Let $\cE_C$ denote the space of pairs $(L,\Phi)$, where $L \in \Pic^1(C)$ is a line bundle and $\Phi: \bP(h_*L) \to \bF_g$ is an orientation-preserving diffeomorphism (not necessarily holomorphic). There is a distinguished component $\cE_{C,0}$ containing $E_C$, i.e. those pairs $(\cO(p), \Phi)$, where $\Phi: \bP(h_*\cO(p)) \to \bF_g$ is a biholomorphism, as in \eqref{eq:EC}. 
\end{definition}

\begin{lemma}\label{lemma:epsextend}
    There is a continuous map
    \[
    \tilde \epsilon: \cE_C \to \Emb(C,\bF_g)
    \]
    extending the embedded monodromy $\tilde \epsilon: E_C \to \Emb(C,\bF_g)$.
\end{lemma}
\begin{proof}
    By \Cref{lemma:Mexists}, the projection $\cE_C \to \Pic^1(C)$ pulls back the tautological $\bF_g$-bundle over $\Pic^1(C)$ to $\cE_C$; call this bundle $\sU$. The data encoded by $\Phi: \bP(h_*L) \to \bF_g$ yields a smooth trivialization 
    \[
    \sU \cong \cE_C \times \bF_g.
    \]
    By construction, there is a fiberwise embedding $C \times \cE_C \to \sU$; composing this with the projection $\sU \cong \cE_C \times \bF_g \to \bF_g$ yields the claimed map.
\end{proof}

By construction, $\cE_C$ is a principal $\Diff^+(\bF_g)$-bundle over $\Pic^1(C)$, and $\cE_{C,0}$ is a principal $G$-bundle over $\Pic^1(C)$, for $G \le \Diff^+(\bF_g)$ some subgroup consisting of a union of connected components. We have the following description of the fundamental group.

\begin{lemma}\label{lemma:centralext}
    There is a central extension
    \[
    1 \to \pi_1(\Diff_0(\bF_g)) \to \pi_1(\cE_{C,0}) \to H_1(C;\Z) \to 0,
    \]
    where as usual, $\Diff_0(\bF_g) \le \Diff(\bF_g)$ denotes the identity component.
\end{lemma}
\begin{proof}
    Taking basepoint at $\id \in \Diff_0(\bF_g) \le G$, the fiber sequence for $G\to \cE_{C,0} \to \Pic^1(C)$ contains the segment
    \[
    \pi_2(\Pic^1(C)) \to \pi_1(\Diff_0(\bF_g)) \to \pi_1(\cE_{C,0}) \to \pi_1(\Pic^1(C)) \to \pi_0(G).
    \]
    The leftmost term vanishes, and the map $\pi_1(\cE_{C,0}) \to \pi_1(\Pic^1(C))$ is a surjection, since the surjection $\pi_1(E_C) \to \pi_1(C) \to \pi_1(\Pic^1(C))$ factors through it (the latter map $\pi_1(C) \to \pi_1(\Pic^1(C))$ is the Abel-Jacobi map). Thus the above sequence is short exact as claimed, and we moreover deduce $G = \Diff_0(\bF_g)$. The extension is central, as this is always the case for a principal $G$-bundle over a path-connected base whenever $G$ is path-connected (recall \cite[Fact 1.4A]{nori}). 
\end{proof}

\subsection{Embedded monodromy detects algebraic intersection} Following \Cref{lemma:epsextend} and \Cref{lemma:centralext}, we see that $\epsilon_*: \pi_1(U_{g,g}) \to \pi_1(\Gr(C,\bF_g))$ cannot be injective, since its restriction to $\pi_1(E_C)$ factors through the group $\pi_1(\cE_{C,0})$, which is two-step nilpotent by \Cref{lemma:centralext}. Thus the best that can be hoped for is for $\tilde \epsilon_*$ to detect elements in the commutator subgroup of $\pi_1(E_C)$. The next result \Cref{lemma:inontriv} shows that $\tilde \epsilon_*$ does detect certain elements in this group, while the subsequent \Cref{lemma:itriv} shows that it does not detect others.

\begin{lemma}\label{lemma:inontriv}
    Let $\alpha, \beta \in \pi_1(E_C)$ be given, and let $\bar \alpha, \bar \beta$ denote their images under the projection $\pi_1(E_C) \to \pi_1(C)$. If the algebraic intersection $\pair{\bar \alpha,\bar \beta}$ is nonzero, then $\tilde \epsilon_*([\alpha, \beta])$ has infinite order in $\pi_1(\Emb(C, \bF_g))$.
\end{lemma}
\begin{proof} We divide the proof into four steps.

    \para{Step 1: Setup} Suppose towards a contradiction that $\tilde \epsilon_*([\alpha,\beta]) = 0$. Let $i_\alpha: C \times S^1 \to \bF_g \times S^1$ and $i_\beta: C \times S^1 \to \bF_g \times S^1$ be the fiberwise embeddings. If $\tilde \epsilon_*(\alpha)$ and $\tilde \epsilon_*(\beta)$ commute as elements of $\pi_1(\Emb(C, \bF_g))$, then this is witnessed by the existence of a map $I: T^2 \to \Emb(C, \bF_g)$ that restricts to $i_\alpha$ on $S^1 \times \{0\}$ and to $i_\beta$ on $\{0\} \times S^1$, yielding a fiberwise embedding
    \[
    I: C \times T^2 \to \bF_g \times T^2.
    \]
    More generally, if $\tilde \epsilon_*([\alpha, \beta]^k) = 0$ for some $k \ge 1$, then this is witnessed by the existence of a map
    \[
    I: (T^2)^{\# k} \to \Emb(C,\bF_g)
    \]
    such that on each $T^2$ factor, the longitude $\ell_i$ maps via $i_\alpha$ and the meridian $m_i$ maps by $i_\beta$. Set $\Sigma = (T^2)^{\# k}$.
    We will arrive at a contradiction by examining the Poincar\'e dual class
    \[
    \delta \coloneq PD(I(C \times \Sigma)) \in H^2(\bF_g \times \Sigma).
    \]

    \para{Step 2: Computing the Poincar\'e dual}  Assume rational coefficients and write
    \[
    H^2(\bF_g) = \Q\pair{x,y},
    \]
    where $y = PD(F)$ is Poincar\'e dual to the class of the ruling fiber, and $x = s + \tfrac{g}{2}y$, where $s = PD(S)$ is dual to the negative section. These satisfy $x^2 = y^2 = 0$. Then $H^2(\bF_g \times \Sigma) = \Q\pair{x,y,\tau}$, where $\tau$ is pulled back from a generator of $H^2(\Sigma)$, normalized so that $\pair{xy\tau, [\bF_g \times \Sigma]} = 1$.

    Write $\delta = a x + b y + c \tau$ for $a,b,c \in \Q$.
    In light of the relations $S^2 = -g, F^2 = 0, SF =1$ and the definition of $x$, since $C \in \abs{2S+(2g+1)F}$, we compute $a = 2$ and $b = g+1$:
    \[
    \delta = 2x + (g+1)y + c \tau.
    \]
    To determine $c$, we make use of the relation $I^*(\delta) = e(\nu)$, where $\nu$ is the normal bundle of $I(C \times \Sigma)$ in $\bF_g \times \Sigma$. The decomposition
    \[
    I^*T(\bF_g \times \Sigma) \cong T(C \times \Sigma) \oplus \nu
    \]
    induces the relation
    \[
    I^*(p_1(T(\bF_g \times \Sigma))) = p_1(T(C \times \Sigma)) + p_1(\nu).
    \]
    As $\bF_g$ has zero signature, $p_1(T\bF_g) = 0$, and $p_1(T\Sigma) = 0$ for dimension reasons. Thus $p_1(T(\bF_g \times \Sigma)) = 0$, and likewise $p_1(T(C \times \Sigma)) = 0$. Meanwhile, $p_1(\nu) = e(\nu)^2$, from which we deduce 
    \[
    e(\nu)^2 = I^*(\delta)^2 = 0.
    \]

    Applying the Poincar\'e duality relation $\pair{\delta \omega, [\bF_g \times \Sigma] } = \pair{I^*(\omega),[C \times \Sigma]}$ to the class $\omega = \delta^2 \in H^4(\bF_g \times \Sigma)$ shows that $\delta^3 = 0$. As $\delta^3 = 12(g+1)c\, xy\tau$, it follows that $c = 0$, so that
    \[
    \delta = 2x + (g+1)y.
    \]
    \para{Step 3: Pulling back} We now compute $I^*(\delta) \in H^2(C \times \Sigma)$. Evidently we must determine $I^*(x)$ and $I^*(y)$. 
    To compute $I^*(x)$, write 
    \[
    I^*(x) = a \eta + \lambda + b \tau,
    \]
    for $a,b \in \Q$ and $\lambda \in H^1(C) \otimes H^1(\Sigma)$. To determine $\lambda$, write
    \[
    \lambda = \sum_{i = 1}^k z_i \otimes PD(\ell_i) + w_i \otimes PD(m_i),
    \]
    where $z_i, w_i \in H^1(C)$ are to be determined. To determine $z_i$, restrict $I: C \times \Sigma \to \bF_g \times \Sigma$ to $C \times m_i \subset C \times \Sigma$, so that $PD(\ell_i)$ restricts to the generator of $H^1(m_i)$. Over $m_i$, the fiberwise embedding $C \times S^1 \to \bF_g \times S^1$ is classified by $\beta \in \pi_1(E_C)$, so that (with $\sM$ as in \Cref{lemma:Mcohomology}) there is a diagram of fiberwise embeddings
    \begin{equation}\label{eq:diamond}
    \xymatrix{
    & C \times E_C \ar[dr]^{\sI'} & \\
    C \times S^1 \ar[ur]^{\id \times \beta} \ar[dr]_{i_\beta} & & \sM.\\
    & \bF_g \times S^1 \ar[ur]_{j_\beta} &
    }
    \end{equation}
    Here, $\sI': C \times E_C \to \sM$ is the composition of $\sI: C \times \Pic^1(C) \to \sM$ with the map $C \times E_C \to C \times \Pic^1(C)$ given on the second factor by the projection and inclusion $E_C \to F_C \subset \Pic^1(C)$, and $j_\beta: \bF_g \times S^1 \to \sM$ is the map of $\bF_g$-bundles classified by $\beta$. 

    Recall from \Cref{lemma:Mcohomology} that $\xi \in H^2(\sM)$ is the Chern class of $\cO_\sM(1)$, and that $\gamma \in H^1(C)\otimes H^1(\Pic^1(C))$ corresponds to the identity under the identification $H^1(C)\otimes H^1(\Pic^1(C)) \cong \End(H^1(C))$. Pulling $\xi$ back along the top path in the diagram then yields
    \[
    (\id \times \beta)^* \sI'^*(\xi) = \eta + (\id \times \beta)^*(\gamma) = \eta + PD(\bar \beta) \otimes \sigma,
    \]
    where $\sigma \in H^1(S^1)$ denotes a generator.
    On the other hand, $H^2(\bF_g \times S^1) \cong H^2(\bF_g)$, and by \Cref{lemma:xiPDS}, the pullback of $\xi$ to $H^2(\bF_g)$ is $PD(S)$. Thus the pullback along the bottom is given as 
    \[
    i_\beta^* j_\beta^*(\xi) = i_\beta^*(x - \tfrac{g}{2}y) = i_\beta^*(x) - g \eta.
    \]
    Thus the component in $H^1(C) \otimes H^1(S^1)$ is the pullback of $\lambda$ along the inclusion $i_{m_i}: C \times S^1 \to C \times \Sigma$ induced by the inclusion $m_i \subset \Sigma$. This shows that $z_i = PD(\bar \beta)$; an analogous argument shows that $w_i = -PD(\bar \alpha)$ (relative to orientations on $\ell_i, m_i$ such that $\pair{\ell_i, m_i} = 1$):
    \begin{equation}\label{eq:lambda}
         \lambda = \sum_{i = 1}^k PD(\bar \beta) \otimes PD(\ell_i) - PD(\bar \alpha) \otimes PD(m_i).
    \end{equation}
   
    The computation of $I^*(y)$ is similar. Write
    \[
    I^*(y) = 2 \eta + \mu + c \tau
    \]
    for $c \in \Z$ and $\mu \in H^1(C) \otimes H^1(\Sigma)$. By \Cref{lemma:Mcohomology}, pulling $y \in H^2(\sM)$ back along the top row of \eqref{eq:diamond} yields
    \[
    (\id \times \beta)^* \sI'^*(y) = 2 \eta, 
    \]
    while pulling back along the bottom yields $i_{m_i}^*I^*(y)$. Running this analysis for each $m_i, \ell_i$ shows that $\mu = 0$. Then the relation $y^2 = 0$ once again pulls back to imply $c = 0$:
    \[
    I^*(y) = 2 \eta.
    \]

    \para{Step 4: Deriving a contradiction}
    As $I^*(\delta^2) = I^*(4(g+1)xy) = 0$, we conclude from the previous step that $8(g+1)b\, \eta \tau = 0$ and hence $b = 0$. Pulling back the relation $x^2 = 0$ then shows that $I^*(x)$ has mixed term $\lambda \in H^1(C) \otimes H^1(\Sigma)$ satisfying
    \begin{equation*}
        \label{eq:musquared}
        \lambda^2 = 0.
    \end{equation*}
    But computing $\lambda^2$ from \eqref{eq:lambda} shows that 
    \[
    \lambda^2 = -2k \pair{\bar \alpha, \bar \beta} \eta \tau,
    \]
    a contradiction.
   \end{proof}

This next lemma shows that $\tilde \epsilon_*$ {\em only} detects the algebraic intersection number of the factors in a commutator.

\begin{lemma}\label{lemma:itriv}
    Let $\alpha, \beta \in \pi_1(E_C)$ be given, and let $\bar \alpha, \bar \beta \in \pi_1(C)$ denote their images under the projection $\pi_1(E_C) \to \pi_1(C)$. If the algebraic intersection $\pair{\bar \alpha,\bar \beta}$ is zero, then $\tilde \epsilon_*([\alpha, \beta]) = 0$ in $\pi_1(\Emb(C, \bF_g))$.
\end{lemma}
\begin{proof}
As $\pi_1(\Pic^1(C)) = H_1(C;\Z)$ is abelian, the images $[\bar \alpha], [\bar \beta]$ of $\alpha, \beta$ in $\pi_1(\Pic^1(C))$ commute. Thus there is a map $J: T^2 \to \Pic^1(C)$ extending the maps $S^1 \to \Pic^1(C)$ defined by $[\bar \alpha]$ and $[\bar \beta]$. We will show that the pullback of the family 
$\sM \to \Pic^1(C)$ along $J: T^2 \to \Pic^1(C)$ is trivial, yielding a product family $\bF_g \times T^2 \to T^2$. The universal fiberwise embedding $\sI: C \times \Pic^1(C) \to \sM$ then restricts to give a fiberwise embedding $C \times T^2 \to \bF_g \times T^2$, showing that the loops of embeddings given by the meridian and longitude of $T^2$ commute.  

Let these loops be denoted by $\alpha',\beta' \in \pi_1(\cE_{C,0})$. By construction, such loops coincide with $[\bar \alpha], [\bar \beta]$ in $\pi_1(\Pic^1(C))$, and so $\alpha' = \alpha\zeta$ as elements of $\pi_1(\cE_{C,0})$, for some $\zeta \in \pi_1(\Diff_0(\bF_g))$; there is a similar expression relating $\beta$ and $\beta'$. By \Cref{lemma:centralext}, $\pi_1(\Diff_0(\bF_g))$ is central, so that $[\alpha,\beta] = [\alpha',\beta']$. Thus, it will be sufficient to show that the pullback of the family 
$\sM \to \Pic^1(C)$ along $J: T^2 \to \Pic^1(C)$ is trivial.

To show that $J^*(\sM) \cong \bF_g \times T^2$ is trivial, we consider the pullback $(\id \times J)^*(\sV)$, a rank-$2$ complex vector bundle over $\P^1 \times T^2$ (here $\sV = (h \times \id)_*(\sL)$ as in \Cref{lemma:Mexists}). By \Cref{lemma:chV}, the Chern classes can be computed as 
\[
c_1((\id \times J)^*(\sV)) = -gy,\qquad c_2((\id \times J)^*(\sV)) = \pair{\bar \alpha, \bar \beta}y \tau,
\]
where $\tau \in H^2(T^2)$ is a generator. The computation of $c_2$ relies on the formula $J^*(\theta) = \pair{\bar \alpha, \bar \beta} \tau$, where $\theta \in H^2(\Pic^1(C))$ is the class of the polarization, represented explicitly under the identification $H^2(\Pic^1(C)) \cong \wedge^2 H_1(C)$ as
\[
\theta = \sum_{i = 1}^g x_i \wedge y_i,
\]
where $\{x_1, y_1, \dots, x_g, y_g\}$ is a symplectic basis for $H_1(C)$. 

Consequently, under the assumption $\pair{\bar \alpha, \bar \beta} = 0$, we conclude that $c_2((\id \times J)^*(\sV)) = 0$. As $c_2$ equals the Euler class, this shows the primary obstruction for $(\id \times J)^*(\sV)$ to admit a smooth nonvanishing section vanishes, and as the base $\P^1 \times T^2$ has dimension four, all higher obstruction classes vanish; we conclude that there is a smooth splitting
\[
(\id \times J)^*(\sV) \cong \underline{\C} \oplus L
\]
for some complex line bundle $L$ on $\P^1 \times T^2$. As these are classified topologically by their Chern class, and $c_1((\id \times J)^*(\sV)) = -gy$, we conclude that $L \cong p^*(\cO_{\P^1}(-g))$, where $p: \P^1 \times T^2 \to \P^1$ is the projection. 

Altogether then, we conclude 
\[
(\id \times J)^*(\sV) \cong p^*(\underline{\C} \oplus \cO_{\P^1}(-g)).
\]
Projectivizing, we find that the pullback of $\sM \to \Pic^1(C)$ along $J: T^2 \to \Pic^1(C)$ is trivial, as desired.
\end{proof}

\subsection{Concluding the proof} We now can give a complete description of the image of $\tilde \epsilon_*$ on the subgroup $\pi_1(E_C) \le \pi_1(U_{g,g})$, and conclude the proof of \Cref{main:Fg} in its more precise form \Cref{thm:preciseFg}.

To proceed, recall that the {\em lower central series} of a group $G$ is defined inductively by $G_{(1)} = G$ and $G_{(k+1)} = [G,G_{(k)}]$. By definition, $G$ is $k$-step nilpotent if $G_{(k+1)} = \{1\}$ is trivial. Note also that each successive quotient $G_{(k)}/G_{(k+1)}$ is abelian.

\begin{lemma}\label{lemma:2mod3}
    There is an isomorphism
    \[
    \pi_1(E_C)_{(2)}/\pi_1(E_C)_{(3)} \cong \wedge^2 H_1(C;\Z).
    \]
\end{lemma}
\begin{proof}
    Recall from \Cref{lemma:ECg} that $\pi_1(E_C)$ is given as a central extension
    \[
    1 \to A \to \pi_1(E_C) \to \pi_1(C) \to 1,
    \]
    with $A = \Z\pair{v,w}/(v^g = w^2)$ and with Euler class $v^2$. Thus $\pi_1(E_C)$ admits the presentation
    \[
    \pi_1(E_C) \cong \pair{x_1, y_1, \dots, x_g, y_g, v,w \mid \prod_{i = 1}^g [x_i,y_i]= v^2,\ v^g = w^2,\ v,w\mbox{ central}}.
    \]
    Since $v,w$ are central, $\pi_1(E_C)_{(2)}/\pi_1(E_C)_{(3)}$ is an abelian group generated by the commutators $[z_p, z_q]$ for $z_p,z_q \in \{x_1,y_1 \dots, x_g,y_g\}$, and the presentation for $\pi_1(E_C)$ shows that no further relations are imposed, so that the map $\wedge^2 H_1(C;\Z) \to \pi_1(E_C)_{(2)}/\pi_1(E_C)_{(3)}$ given by $z_p \wedge z_q \mapsto [z_p,z_q]$ is an isomorphism.
\end{proof}

\begin{theorem}\label{thm:preciseFg}
    The embedded monodromy $\tilde \epsilon_*: \pi_1(E_C) \to \pi_1(\Emb(C, \bF_g))$ factors through the quotient $\pi_1(E_C)/\pi_1(E_C)_{(3)}$ by the third term in the lower central series. Under the isomorphism
    $\pi_1(E_C)_{(2)}/\pi_1(E_C)_{(3)} \cong \wedge^2 H_1(C;\Z)$ of \Cref{lemma:2mod3}, the restriction of $\tilde \epsilon_*$ can be identified with the contraction homomorphism
    \begin{align*}
     \kappa: \wedge^2 H_1(C;\Z) &\to \Z\\
     x \wedge y &\mapsto \pair{x,y}.
    \end{align*}
The image $\tilde \epsilon_*(\pi_1(E_C)) \le \pi_1(\Emb(C, \bF_g))$ is a central extension of $H_1(C;\Z)$ by an abelian group of rank one. 
\end{theorem}

\begin{proof}
    By \Cref{lemma:epsextend}, the map $\tilde \epsilon: E_C \to \Emb(C,\bF_g)$ extends to $\tilde \epsilon: \cE_{C,0} \to \Emb(C, \bF_g)$. At the level of $\pi_1$, this shows that $\tilde\epsilon_*: \pi_1(E_C) \to \pi_1(\Emb(C,\bF_g))$ factors through $\pi_1(\cE_{C,0})$. By \Cref{lemma:centralext}, this latter group is two-step nilpotent, showing that $\tilde\epsilon_*$ factors through $\pi_1(E_C)/ \pi_1(E_C)_{(3)}$ as claimed.

    The image of $\tilde\epsilon_*$ on the abelian group $\wedge^2 H_1(C;\Z)$ is determined by its image on a generating set, of the form $z_p \wedge z_q$ for $z_p,z_q \in \{x_1, y_1, \dots, x_g, y_g\}$ as in \Cref{lemma:2mod3}. By \Cref{lemma:itriv}, it vanishes on any $z_p \wedge z_q$ with $\pair{z_p,z_q} = 0$, leaving us to determine only the values $\tilde\epsilon_*(x_i \wedge y_i)$. But as $\pair{x_i + x_j, y_i - y_j} = 0$, it follows that $
    \tilde\epsilon_*(x_i \wedge y_i) = \tilde \epsilon_*(x_j \wedge y_j)$ for all pairs of indices $1 \le i,j\le g$, and by \Cref{lemma:inontriv}, this common element has infinite order; we denote it by $\zeta$.

    To fully describe the image of $\tilde \epsilon_*$, we let $q: \pi_1(E_C) \to H_1(C;\Z)$ be the natural quotient. The work of the previous paragraph establishes the formula
    \[
    [\tilde \epsilon_*(\alpha), \tilde \epsilon_*(\beta)] = \zeta^{\pair{q(\alpha),q(\beta)}},
    \]
    valid for all $\alpha, \beta \in \pi_1(E_C)$. As the intersection pairing is nondegenerate and $\zeta$ is of infinite order, it follows that $\ker(\tilde \epsilon_*) \le \ker(q)$, so that $q$ factors through $\tilde \epsilon_*$ and hence there is an extension
    \[
    1 \to K \to \tilde \epsilon_*(\pi_1(E_C)) \to H \to 1.
    \]
    for some subgroup $K$. These same observations in fact show that $K = Z(\tilde \epsilon_*(\pi_1(E_C)))$ is the center: an element $\tilde \epsilon_*(\alpha)$ is central exactly when the functional $\pair{q(\alpha), \cdot}$ is zero. The presentation of $\pi_1(E_C)$ obtained in \Cref{lemma:Apres} then shows that $K$ is generated by the elements $\zeta, v,w$, subject to the relations $g \zeta = 2v$ and $gv = 2w$ (written additively). Possibly there are other relations as well, but by \Cref{lemma:inontriv}, $\zeta \in K$ has infinite order, and so whatever additional relations may hold, $K$ must be an abelian group of rank one, as claimed.
\end{proof}

\begin{proof}[Proof of \Cref{cor:diff}]
    Combining \Cref{lemma:ECg} and \Cref{lemma:centralext}, there is a commutative diagram with exact rows
    \[
    \xymatrix{
    1 \ar[r] & A \ar[r] \ar[d] & \pi_1(E_C) \ar[r] \ar[d] & \pi_1(C) \ar[r] \ar[d]\ar[r] & 1\\
    1 \ar[r] & \pi_1(\Diff_0(\bF_g)) \ar[r] & \pi_1(\cE_{C,0}) \ar[r] & H_1(C;\Z) \ar[r] & 1,
    }
    \]
    where the map $\pi_1(C) \to H_1(C;\Z)$ is the abelianization, induced by the Abel-Jacobi map $C \to \Pic^1(C)$. Moreover, $\tilde\epsilon_*: \pi_1(E_C) \to \pi_1(\Emb(C,\bF_g))$ factors through $\pi_1(\cE_{C,0})$. Commutativity and exactness of the above diagram then shows that the restriction of $\tilde\epsilon_*$ to $\pi_1(E_C)_{(2)}$ factors through $\pi_1(\Diff_0(\bF_g))$. By \Cref{thm:preciseFg}, the restriction of $\tilde\epsilon_*$ to $\pi_1(E_C)_{(2)}$ has infinite cyclic image, proving the claim.
\end{proof}

\begin{remark}
In the case of $\bF_0 \cong S^2 \times S^2$, there is a second source of infinite-order elements in $\pi_1(\Diff(\bF_0))$, coming from \Cref{main:P1xP1}. This result shows that the theta meridian element $z \in \pi_1(U_{0,g})$ satisfies $\bar \epsilon_*(z) = 1$ in $\pi_0(\Diff(\bF_0,C))$, but $\tilde \epsilon_*(z) \ne 1$ in $\pi_1(\Emb(C,\bF_0))$. As the intrinsic monodromy $\rho_*(z) \in \Mod(C)$ is trivial, the target $\bar \epsilon_*(z)$ can be lifted from $\pi_0(\Diff(\bF_0, C))$ to the ``pure subgroup'' $\pi_0(\PDiff(\bF_0,C))$ (the component group of diffeomorphisms of $\bF_0$ that fix $C$ pointwise). Under the long exact sequence of homotopy groups for the bundle $\Diff(\bF_0) / \Emb(C, \bF_0)$ (with fiber $\PDiff(\bF_0,C)$), this shows that $\tilde\epsilon_*(z)$ lifts to an infinite-order element of $\pi_1(\Diff(\bF_0))$. 

This can be understood somewhat more explicitly as follows. The theta meridian $z$ induces a loop of embeddings of a fixed hyperelliptic curve $C$ of genus $g$ into $\bF_0$ that are ``close to degenerating'' into an embedding of $C$ into the diffeomorphic algebraic surface $\bF_2$. The loop of diffeomorphisms is constructed by pushing this loop of hyperelliptic curves in the ambient $S^2 \times S^2$, obtaining a diffeomorphism of the pair $(S^2 \times S^2, C)$ by isotopy extension. As the induced relative mapping class in $\pi_0(\Diff(S^2 \times S^2, C))$ is trivial, one can compose this path of diffeomorphisms with an isotopy of the pair $(S^2\times S^2, C)$ back to the identity, completing the loop in $\Diff(S^2 \times S^2)$. To understand this more concretely, one should develop a local model for what is happening to the curves as the divisor parameterizing them limits onto the theta divisor in $\Pic^{g+1}(C)$.
\end{remark}

\bibliography{bib}{}
\bibliographystyle{alpha}
\end{document}